\documentclass[11pt, reqno, a4paper]{amsart}

\usepackage{amsmath, amssymb, amsthm}
\usepackage{amscd}
\usepackage[figuresright]{rotating}
\usepackage{mathtools} %to use mathclap

\usepackage{color}

\usepackage{dsfont}
\usepackage[matrix, arrow, curve, frame, arc]{xy}
\usepackage{pgf}
\usepackage{tikz-cd}
\usetikzlibrary{arrows,shapes,automata,backgrounds,petri,calc}
\usepackage[square,sort,comma,numbers]{natbib}
 \usepackage{url}
 \usepackage{hyperref}
 \usepackage[shortlabels]{enumitem} %to index the enumeration

\newcommand{\tikzAngleOfLine}{\tikz@AngleOfLine}
\def\tikz@AngleOfLine(#1)(#2)#3{%
\pgfmathanglebetweenpoints{%
\pgfpointanchor{#1}{center}}{%
\pgfpointanchor{#2}{center}}
\pgfmathsetmacro{#3}{\pgfmathresult}%
}

\newcommand{\End}{\operatorname{End}}
\newcommand{\Hom}{\operatorname{Hom}}
\newcommand{\gldim}{\operatorname{gldim}}

\newcommand{\Ext}{\operatorname{Ext}}
\newcommand{\Tor}{\operatorname{Tor}}
\newcommand{\add}{\!\operatorname{add}}
\newcommand{\pdim}{\operatorname{pdim}}
\newcommand{\idim}{\operatorname{idim}}
\newcommand{\m}{\!\operatorname{-mod}} 
 
\newcommand{\rmod}{\operatorname{mod-}\!\!}
\newcommand{\proj}{\!\operatorname{-proj}}

\newcommand{\id}{\operatorname{id}}

\renewcommand{\top}{\operatorname{top}}

\newcommand{\im}{\!\operatorname{im}}

\newcommand{\injdim}{\operatorname{idim}}
\newcommand{\domdim}{\operatorname{domdim}}
\newcommand{\R}{\operatorname{R}}

\newcommand{\codomdim}{\operatorname{codomdim}}
 
\newcommand{\cograde}{\operatorname{cograde}}

\newcommand{\Tr}{\operatorname{Tr}} 
\newcommand{\dd}{\!\operatorname{-domdim}}
\newcommand{\cdd}{\!\operatorname{-codomdim}}
\newcommand{\grade}{\operatorname{grade}}
\newcommand{\Dom}{\operatorname{Dom}}
\newcommand{\TF}{\operatorname{TF}}

\newcommand{\otimesL}{\otimes^{\operatorname{\bf L}}}
\newcommand{\RHom}{{\rm\bf R}{\operatorname{Hom}}}

\newtheorem{numberingthm}{Theorem}[section] % reset theorem numbering for each section
\theoremstyle{definition}
\newtheorem{Def}[numberingthm]{Definition}
\newtheorem{example}[numberingthm]{Example}
\newtheorem{conjecture}[numberingthm]{Conjecture}
\newtheorem{question}[numberingthm]{Question}

\theoremstyle{plain}
\newtheorem{Prop}[numberingthm]{Proposition}
\newtheorem{Theorem}[numberingthm]{Theorem}
\newtheorem{Conjecture}[numberingthm]{Conjecture}
\newtheorem{Cor}[numberingthm]{Corollary}
\newtheorem{Lemma}[numberingthm]{Lemma}
\newtheorem{Question}[numberingthm]{Question}
\newtheorem{Remark}[numberingthm]{Remark}
\newtheorem{thmintroduction}{Theorem}

\newtheorem*{nonumbercor}{Corollary}
\newtheorem*{nonumberconj}{Conjecture}

\newenvironment{Example}
{\pushQED{\qed}\example}
{\popQED\endexample}

\theoremstyle{remark}

\begin{document}

\title[]{A new formula for the classical dominant dimension using bimodules}

\author[T. Cruz]{Tiago Cruz}
\address[Tiago Cruz]{Institut f\"ur Algebra und Zahlentheorie,
Universit\"at Stuttgart, Germany }
\email{tiago.cruz@mathematik.uni-stuttgart.de}

\author[R. Marczinzik]{Ren\'{e} Marczinzik}
\address[Ren\'{e} Marczinzik]{Mathematical Institute of the University of Bonn, Endenicher Allee 60, 53115 Bonn, Germany}
\email{marczire@math.uni-bonn.de}
%\date{\today}

\begin{abstract}
We show that a faithful projective-injective module over a finite-dimensional algebra $A$ has the double centraliser property if and only if $A$ as a bimodule is reflexive. More generally, we provide a new characterisation of the classical dominant dimension by showing that having dominant dimension at least $n$ is equivalent to the bimodule $A$ being $n$-torsion-free.
This allows us to find new connections between the classical Tachikawa and Nakayama conjectures and Gorenstein homological algebra.
Furthermore, we use our results to give new interpretations of Hochschild (co)homology of finite-dimensional algebras using higher Auslander-Reiten translates and the canonical bimodule in the sense of Fang, Kerner and Yamagata.
\end{abstract}

\maketitle

\section{Introduction}
Many fundamental dualities in mathematics arise from reflexivity and double centraliser properties and their analogues, where a structure is naturally recovered through a bidual or a double centraliser. Prominent examples appear in linear algebra, functional analysis (for instance, reflexive Banach spaces), ring theory and representation theory (for instance, Morita duality, Schur-Weyl duality, Matlis duality, reflexive modules), in algebraic geometry (for instance, reflexive sheaves), among many others.
Let $A$ be a finite-dimensional $K$-algebra over a field $K$ with right $A$-module $M$.
Then $M$ is a left $B$-module, where $B = \End_A(M)$ and $M$ is said to have the \emph{double centraliser property} in case the canonical homomorphism of algebras $f : A \rightarrow \End_{B^{op}} (M)^{op}$ is an isomorphism, where
$f(a)(x) = xa$. A natural and insightful instance occurs when $M = eA$ for an idempotent $e$ such that $eA$ is a faithful projective-injective $A$-module.
Double centraliser properties, particularly on faithful projective-injective modules, play an important role in the representation theory of algebras and algebraic Lie theory. Notable cases are Soergel's Struktursatz \cite{Soe}, the duality between the representation theory of $q$-Schur algebras $S(n,r)$ for $n \geq r$ and the representation theory of Iwahori-Hecke algebras (see \cite{KSX}, \cite{Do} and the references therein), the higher Auslander correspondence that relates cluster tilting subcategories with higher Auslander algebras \cite{Iya} and the many specialisations of the Morita-Tachikawa correspondence (for instance \cite{zbMATH03248955, FanKoe2, zbMATH06833443}).

Recall that the \emph{dominant dimension} of a module $M$ with minimal injective coresolution
$$0 \rightarrow M \rightarrow I^0 \rightarrow I^1 \rightarrow \cdots $$
is defined as the smallest $n$ such that $I^n$ is not projective and infinite if no such $n$ exists. The dominant dimension of an algebra is then defined as the dominant dimension of the regular module.
The dominant dimension reflects important aspects of the homological structure of an algebra. For example, it can be used to study isomorphisms of Ext-spaces, see for example \cite{FK3}, and it is a defining property in the definition of higher Auslander algebras \cite{Iya}. In some cases, the dominant dimension of modules can be seen as a replacement of the depth for non-commutative finite-dimensional algebras, see \cite{CM25}.
The main goal of this paper is to establish that dominant dimension  encodes bimodule-theoretic properties of $A$, such as reflexivity and its higher analogues, as well as the property of being an $n$-syzygy bimodule. For this goal, the generalisation of dominant dimension studied in \cite{cruz2022split} is the tool connecting all these concepts.
\[
* \ \ * \ \ *
\]

Note that viewing $A$ as an $A$-bimodule is equivalent to viewing $A$ as a left $A^e$-module, where $A^e= A \otimes_K A^{op}$ denotes the enveloping algebra of $A$.
The main result of this article is a new formula for the classical dominant dimension in terms of the relative dominant dimension of the regular bimodule. Namely, if $A^e\dd_{A^e} A$ denotes the relative dominant dimension of $A$ as an $A^e$-module with respect to the regular module $A^e$ (see Subsection  \ref{Torsion-free modules revisited} for the definition), then we prove the following.

\begin{thmintroduction} (see Theorems \ref{mainresult}, \ref{thm2dot3}, \ref{thm2dot7} and Corollary \ref{cornsyzygydominantdimension})\label{mainformulaintro}
Let $A$ be a finite-dimensional algebra. Then
\[
        \domdim A
        =
        A^e\dd_{A^e} A
        =
        \sup\{\, n \in \mathbb{N}_0 \mid A \in \Omega^n(A^e\m)\,\}.
\]
Moreover, for every $n \geq 1$, the following are equivalent:
\begin{enumerate}[(1)]
\item $A$ has dominant dimension at least $n$;
\item $A$ is $n$-torsion-free as a bimodule.
%\item $A^e\dd_{A^e} A \geq n$.
\end{enumerate}
\end{thmintroduction}

Here, $\Omega^n(A^e\text{-mod})$ denotes the category of bimodules arising as $n$-syzygies, possibly up to projective summands, and the terminology of $n$-torsion-free modules is understood in the sense of Auslander--Bridger. Indeed, an $A$-module $M$ is called \emph{reflexive} if the canonical evaluation map $M \rightarrow M^{**}$ is an isomorphism, where $(-)^{*}=\Hom_A(-,A)$. Recall that a famous theorem of Auslander and Bridger \cite{AB69} states that a module $M$ is reflexive if and only if $\Ext_A^i(\Tr M,A)=0$ for $i=1,2$, where $\Tr M$ denotes the Auslander-Bridger transpose of $M$.
This result motivated the generalisation of reflexive modules as follows: $M$ is \emph{$n$-torsion-free} for some $n \geq 1$ if and only if $\Ext_A^i(\Tr M,  A)=0$ for $1 \leq i \leq n.$ We refer for example to \cite{AR3} for more on $n$-torsion-free modules for Artin algebras.

The cases $n=1,2$ in Theorem \ref{mainformulaintro} can be obtained by combining techniques on the canonical bimodule $\Hom_A(D(A), A)$ of Fang, Kerner and Yamagata \cite{zbMATH06814513} with the language of relative dominant dimension. The remaining cases, namely $n\geq 3$, require spectral sequences and the framework of derived categories.
\[
* \ \ * \ \ *
\]
As a first consequence, we recover and extend the connection between double centraliser properties and reflexivity. 
The study of reflexive modules is a classical topic in ring theory and algebraic geometry; we refer, for example, to \cite{Es} and \cite{EK2}.
By the celebrated Morita--Tachikawa correspondence, an algebra $A$ has dominant dimension at least $2$ if and only if it admits a faithful projective-injective module satisfying the double centraliser property, see for example \cite{Yam}. Combining this with Theorem \ref{mainformulaintro} gives the following particular case.

\begin{nonumbercor}(see Corollary \ref{cor3dot7} and Theorem \ref{thm2dot2}) \label{maincorollaryintro}
	Let $A$ be a finite-dimensional algebra. Then the following are equivalent:
\begin{enumerate}[(1)]
\item $A$ has a faithful projective-injective $A$-module which has the double centraliser property.
	\item $A$ as a bimodule is reflexive.
\end{enumerate}
\end{nonumbercor}

The first connection between the reflexivity of modules and their dominant dimension was discovered by Morita in \cite{zbMATH03271644}. This understanding was significantly extended by the works \cite{MV2}, \cite{AR3}, and more recently \cite{zbMATH06814513}. Theorem \ref{mainformulaintro} provides a uniform extension of this connection to all degrees.

As a consequence of our results, we can connect two of the most important homological conjectures for finite-dimensional algebras to Gorenstein homological algebra.
The first conjecture is the \emph{Nakayama conjecture}, which states that a finite-dimensional algebra has infinite dominant dimension if and only if it is self-injective.
The second is the \emph{first Tachikawa conjecture}, which states that an algebra $A$ is self-injective if and only if $\Ext_A^i(D(A),A)=0$ for all $i>0$; see \cite{Tachikawabook}.
This terminology distinguishes it from the \emph{second Tachikawa conjecture}, stating that every self-orthogonal module over a self-injective algebra is projective. The Nakayama conjecture and the Tachikawa conjectures have been studied extensively, we refer for example to the recent articles \cite{MR5047073,MR4956374,MR4843724,MarczinzikOwens2026}.
The Nakayama conjecture and the Tachikawa conjectures are wide open and only known for special classes of algebras, such as monomial algebras, and both are consequences of the finitistic dimension conjecture. We refer for example to \cite{Yam} for the relation of those conjectures to the other classical homological conjecture for finite dimensional algebras.

Gorenstein homological algebra can be seen as a generalisation of classical homological algebra, where instead of projective modules, the Gorenstein-projective modules play the main role.
We refer for example to the survey \cite{Chen}.
A module $M$ is called \emph{Gorenstein-projective} if $M$ is reflexive and $\Ext_A^i(M, A)=\Ext_A^i(M^*, A)=0$ for all $i>0$.
Our second main result shows that, for the regular bimodule, the Gorenstein-projective condition is equivalent to combining infinite dominant dimension with the vanishing condition appearing in the first Tachikawa conjecture.

\begin{thmintroduction} (see Theorem \ref{thm6dot4})\label{secondmainresultintro}
Let $A$ be a finite-dimensional algebra. Then the following are equivalent:
\begin{enumerate}[(1)]
\item $A$ is Gorenstein-projective as an $A$-bimodule;
\item $\domdim A=+\infty$ and $\Ext_A^i(D(A),A)=0$ for all $i>0$;
\item $\Ext_{A^e}^i(A,A^e)=0$ for all $i>0$ and $A^e\dd_{A^e} A=+\infty$.
\end{enumerate}
\end{thmintroduction}

This theorem motivates the following new homological conjecture.

\begin{nonumberconj}
Let $A$ be a finite-dimensional algebra.
Then $A$ as an $A$-bimodule is Gorenstein-projective if and only if $A$ is self-injective.
\end{nonumberconj}
We call this conjecture the \emph{Gorenstein bimodule conjecture}. Recently, Chen and Xi in \cite{MR5047073} proved that this conjecture is true for virtually Gorenstein algebras.

Our methods give us the following result, proved in Theorem \ref{gorbimodtheorem}.
\begin{thmintroduction}
Let $A$ be a given finite-dimensional algebra. Then the following are equivalent:
\begin{enumerate}
\item $A$ satisfies the Gorenstein bimodule conjecture.
\item $A$ satisfies the Nakayama conjecture or the first Tachikawa conjecture.
\end{enumerate}
\end{thmintroduction}
In particular, the Gorenstein bimodule conjecture is a consequence of the finitistic dimension conjecture.

Our final application gives a new formula to calculate the classical Hochschild homology and cohomology for a finite-dimensional algebra.
For this, we use the \emph{canonical bimodule} $V:=\Hom_A(D(A),A)$ introduced by Fang, Kerner and Yamagata in \cite{zbMATH06814513}.
\begin{thmintroduction} (see Theorem \ref{formulahochschild})
Let $A$ be a finite-dimensional algebra with dominant dimension $n \geq 2$.
\begin{enumerate}
\item We have for the Hochschild homology and $l \geq 1$:
$$HH_l(A) \cong D \Ext_{A^e}^{l+n}(A, \tau_{n-1}(V)).$$
\item We have for the Hochschild cohomology and $l \geq 1$:
$$HH^l(A) \cong \Ext_{A^e}^{l+n}(D(A), \tau_{n-1}(V)).$$
\end{enumerate}
\end{thmintroduction}
Here $\tau_d:=\tau \Omega^{d-1}$ denotes the $d$-th \emph{higher Auslander-Reiten translate} with $\tau$ being the classical Auslander-Reiten translate.
The new formula for the Hochschild cohomology raises the question of what the projective dimension of $D(A)$ as an $A$-bimodule is, as this would allow us to obtain new vanishing results for the Hochschild cohomology.
We will see in Corollaries \ref{finitenessglobalusingregularasbimodule} and \ref{boundsforinjdimasbimodule} that $g \leq \pdim_{A^e}(D(A)) \leq 2 g$ when $g$ denotes the global dimension of $A$ and in particular $\pdim_{A^e}(D(A))$ is infinite when the algebra has infinite global dimension. In Proposition \ref{injdimasbimoduleleftrightsymm}, we see that $\pdim_{A^e}(D(A))= \idim_{A^e} A$.
\[
* \ \ * \ \ *
\]
The paper is organised as follows. Section \ref{Preliminaries} introduces the necessary notation to be used throughout, reductions of cohomology between bimodules to cohomology between one-sided complexes and spectral sequences. Section \ref{On torsion-free modules and Gorenstein projective modules} recalls the relevant material on dominant dimension, relative dominant dimension, $n$-torsion-free modules and Gorenstein-projective modules. In Section \ref{Bimodule characterisations of the dominant dimension}, we prove the main equality of the paper, characterising the classical dominant dimension in terms of the relative dominant dimension of the regular bimodule, and derive the corresponding $n$-torsion-free and syzygy characterisations. In Subsection \ref{Revisiting Mueller's theorem and gendo-symmetric algebras}, we give an alternative proof of Mueller's characterisation of dominant dimension for gendo-symmetric algebras, following the formulation given by Fang and Koenig. In Section \ref{Applications to Hochschild (co)homology}, we present applications of our main result to Hochschild homology and Hochschild cohomology, and we prove that $A$ has finite injective dimension as a bimodule if and only if it has finite global dimension.
Section \ref{Applications related to homological conjectures} deals with applications to the homological conjectures, while in Section \ref{Further questions} we indicate some open problems for general noetherian rings motivated by results in this article.

\section{Preliminaries} \label{Preliminaries}
We will always assume that $A$ is a finite-dimensional $K$-algebra over a field $K$ and modules are finitely generated unless otherwise stated. To distinguish left and right $A$-modules, we write ${}_A M$ for a left $A$-module and $M_A$ for a right $A$-module. When no confusion arises, we omit the side. We denote by $\rmod A$ the category of finitely generated right $A$-modules, and by $A\m$ the category of finitely generated left $A$-modules. We write $D(-)=\Hom_K(-,K)$ to denote the standard duality $\rmod A\longleftrightarrow \rmod A^{op}$, where $A^{op}$ denotes the opposite algebra of $A$. We denote by $\pdim_A M$ and $\injdim_A M$ the projective and injective dimension of an $A$-module $M$, respectively. When there is no ambiguity, we may omit the subscript $A$. The global dimension of the algebra $A$ is denoted by $\gldim A$.
For an introduction to representation theory and homological algebra of finite-dimensional algebras, we refer for example to \cite{ASS}, \cite{SkoYam} and \cite{Z}. For an introduction to Gorenstein homological algebra, we refer for example to the survey \cite{Chen}.

\subsection{Bimodules, Ext and Tor}

%TOCITE: Cartan-Eilenberg, Chapter 9, Corollary 4.4: $\Ext_{A^e}^n(A, \Hom_k(B,C)) \cong \Ext_A^n(B,C)$ and use $\Hom_k(B,C) \cong \Hom_k(B,k) \otimes_k C$
%something introductory to bimodules
In this subsection, we present some homological properties of bimodules and how they descend to single-sided modules.

We write $A^e=A\otimes_K A^{op}={}_AA\otimes_K A_A$ for the \emph{enveloping algebra} of $A$. Over a finite-dimensional algebra $A$-bimodules coincide with the left $A^e$-modules. Recall that a right $A^e$-module $X$ has a left $A$-module structure through the action $a\cdot x:=x (1\otimes a)$. 
Observe that $(A^e)^{op}\cong A^e$ via $x\otimes y\mapsto y\otimes x$. So, from a representation-theoretical point of view, it is equivalent to work with left $A^e$-modules or with right $A^e$-modules.
In the following, we can interchange the results between right and left modules due to the following fact: $M\otimes_{A^{op}} N\cong N\otimes_A M$.

\begin{Prop} \label{torformula}
Let $A$ be a finite-dimensional algebra with a left $A$-module $M$ and a right $A$-module $N$.
\begin{enumerate}
\item $\Tor_i^A(N,M) \cong D \Ext_A^i(M,D(N))$ for all $i \geq 0$.
\item If $P$ is a projective $A$-bimodule, then $P \otimes_A M$ and $N\otimes_A P$ are projective $A$-modules.

\end{enumerate}
\end{Prop}
\begin{proof}
\begin{enumerate}
\item See, for example, Proposition 4.11. in appendix A of \cite{ASS}.
\item See, for example, Lemma 11.15. in chapter IV. in \cite{SkoYam} or \citep[Lemma 2.1]{zbMATH00067188}.
\end{enumerate}
\end{proof}
As we will see below, relating Ext groups between bimodules and one-sided modules involves complexes. Thus, the derived category provides the natural setting for carrying out such comparisons. So, we shall recall some conventions that we use about complexes.

A complex $X$ of $A$-modules is a sequence of morphisms $d^i$ in $A\m$ (or in $\rmod A$) of the form $\cdots\rightarrow X^i\xrightarrow{d^{i-1}}X^i \xrightarrow{d^i} X^{i+1} \rightarrow \cdots$ with $d^i\circ d^{i-1}=0$ for all $i\in \mathbb{Z}$. We denote by $\mathcal{D}(A)$ the derived category of finitely generated $A$-modules. Given a subcategory $\mathcal{C}$ of $A\m$, the homotopy category of complexes over $\mathcal{C}$ is denoted by $\mathcal{K}(\mathcal{C})$. For complexes $X, Y\in \mathcal{D}(A)$, we denote the right derived Hom complex by $\RHom_A(X, Y)$ and the left derived tensor product by $X\otimesL_A Y$. For a complex $X$, the shift of $X$ by an integer $i$ is denoted by $X[i]$ and it is defined by $(X[i])^n=X^{n+i}$ with differentials adjusted by a sign. In particular, $H^i(X)=H^0(X[i])$ for all complexes $X$. The $A$-modules are regarded in $\mathcal{D}(A)$ as the complexes where all the terms except the one at index 0 are zero. Moreover,
for $X, Y\in \mathcal{D}(A)$ and $i\in \mathbb{Z}$, the $i^{th}$-extension group of $X$ by $Y$ is $\Ext_A^i(X, Y)=\Hom_{\mathcal{D}(A)}(X, Y[i]).$ For general details on derived categories we refer to \cite{zbMATH00040982}.

\begin{Prop} \label{restrictingfromenvelopingtoregular}
Let $A$ be a finite-dimensional algebra over a field $K$. 
\begin{enumerate}
\item $A^e\cong \Hom_K(D(A), A)$ as $A^e$-bimodules.
\item Let $X$ be a left $A^e$-module. Then, $\Hom_{A^e}(X, A^e)\cong \Hom_A(D(A)\otimes_A X, A_A)$ as right $A^e$-modules. Here, $DA\otimes_A X$ is a right $A$-module through the right $A$-module structure on $X$ and $\Hom_A(D(A)\otimes_A X, A_A)$ is a bimodule over $A$ via the left $A$-module structures on $D(A)$ and $A$.
\item Let $X$ be a right $A^e$-module. Then, $\Hom_{A^e}(X, A^e)\cong \Hom_A(X\otimes_A D(A), {}_AA)$ as left $A^e$-modules.
\item $\RHom_{A^e}(A, A^e)\cong \RHom_A(D(A), A)$ as objects in $\mathcal{D}(A^e)$ and by restriction as objects in $\mathcal{D}(A)$. In particular, this isomorphism holds when $A$ in the first component is regarded as a right $A^e$-module or as a left $A^e$-module. 
\item $\Ext_{A^e}^i(A, A^e)\cong \Ext_A^i(D(A), A)$ for every $i\geq 0$.
\item Let $X$ be a left $A^e$-module. Then, $\Ext^i_{A^e}(X, A^e)\cong \Ext^i_A(D(A)\otimesL_A X, A_A)$ as right $A^e$-modules.
\item Let $X$ be a right $A^e$-module. Then, $\Ext_{A^e}^i(X, A^e)\cong \Ext^i_A(X\otimesL_A D(A), {}_AA)$ as left $A^e$-modules.
\end{enumerate}
\end{Prop}
\begin{proof}
(1) follows from \citep[Corollary 4.4]{zbMATH00067188}. (2) and (3) follow from \citep[Corollary 4.2]{zbMATH00067188}. We now prove (4). Let $A^\bullet$ be a projective resolution of the left $A^e$-module $A$. Then, by (2) and Proposition \ref{torformula},
\begin{align}
\RHom_{A^e}(A, A^e)=\Hom_{A^e}(A^\bullet, A^e)\cong \Hom_A(D(A)\otimes_A A^\bullet, A)= \RHom_A(D(A), A_A)
\end{align} since $\Tor_i^A(D(A), A)=0$ for $i>0$ and $D(A)\otimes_A A^\bullet$ is a projective resolution of the left $A$-module $D(A)$. Analogously, replacing $A$ with $A^{op}$, the result holds when $A$ is regarded as a right $A^e$-module.
The cohomology of the objects in (4) gives (5). To establish (6), let $X$ be a left $A^e$-module and $P^{\bullet}$ be a projective $A^e$-resolution of $X$. Then $DA\otimes_A P^{\bullet}$ is a complex of projective $A$-modules and we have the following
\begin{align*}
	\Ext_{A^e}^i(X, A^e)&\cong \Hom_{\mathcal{D}(A^e)}(X, A^e[i])\cong \Hom_{\mathcal{D}(A^e)}(X[-i], A^e)\cong \Hom_{\mathcal{D}(A^e)}(P^{\bullet}[-i], A^e)\\&\cong \Hom_{\mathcal{K}(A^e\proj)}(P^{\bullet}[-i], A^e)\cong \Hom_{\mathcal{K}(A\proj)}(DA\otimes_A P^{\bullet}[-i], A_A)\\&\cong \Hom_{{\mathcal{D}(A)}}(DA\otimes_A P^{\bullet}[-i], A_A) \cong \Hom_{{\mathcal{D}(A)}}(DA\otimesL_A X[-i], A_A)\\&\cong \Ext^i_A(D(A)\otimesL_A X, A_A).
\end{align*}

%The following identifications follow from the derived Tensor-Hom adjunction:
%\begin{align*}
%\Ext_{A^e}^i(X, A^e)&\cong \Hom_{\mathcal{D}(A^e)}(X, A^e[i])\cong \Hom_{\mathcal{D}(A^e)}(X[-i], A^e) \cong \Hom_{\mathcal{D}(A^e)}(X[-i]\otimesL_A A, A^e)
%\\ &\cong \Hom_{\mathcal{D}(A)}(X[-i], \RHom_{A^e}(A, A^e)) \cong \Hom_{\mathcal{D}(A)}(X[-i], \RHom_{A}(D(A), A_A)) \\&\cong \Hom_{\mathcal{D}(A)}(D(A)\otimesL_A X[-i], A_A)\cong \Hom_{\mathcal{D}(A)}(D(A)\otimesL_A X, A[i]).
%\end{align*}
Applying (6) to the opposite algebra $A^{op}$, we obtain (7).
\end{proof} 

A particularly important case of the previous proposition is obtained by applying it to the dual bimodule $A^*:=\Hom_{A^e}(A, A^e)$, which coincides with the canonical bimodule studied in  \cite{zbMATH06814513}.

\begin{Lemma}\label{lem3dot8}
For every $i\in \mathbb{N}_0$ we have $\Ext_{A^e}^i(A^*, A^e)\cong \Ext_A^i(D(A)\otimesL_A A^*, A)$, where $A^*=\Hom_{A^e}(A, A^e)\in A^e\m$.
In particular, $$\Hom_{A^e}(A^*, A^e)\cong \Hom_A(D(A)_A, \Hom_A(\Hom_A(D(A), A), A_A)).$$
\end{Lemma}
\begin{proof} Applying $i=0$ to Proposition \ref{restrictingfromenvelopingtoregular}(5) yields that $A^*:=\Hom_{A^e}(A, A^e)$ is isomorphic to $\Hom_A(D(A), A)$ as $A$-bimodules. So, the claim follows from Proposition \ref{restrictingfromenvelopingtoregular} (6) with $X=A^*$.

The case $i=0$ yields, in particular, 
\begin{align}
\Hom_{A^e}(A^*, A^e)&\cong \Hom_{\mathcal{D}(A)}(D(A)\otimesL_A A^*, A)\cong \Hom_{\mathcal{D}(A)}(D(A), \RHom_A(A^*, A))\\&\cong \Hom_A(D(A), \Hom_A(A^*, A))\\&\cong \Hom_A(D(A), \Hom_A(\Hom_A(D(A), A), A)).
\end{align}
\end{proof}

Recall that the \emph{grade} of a module $M$ is defined as  $\grade M:=\inf \{ i \geq 0 \mid \Ext_A^i(M,A) \neq 0 \}$.
Dually the \emph{cograde} of a module $M$ is defined as $\cograde M:=\inf \{ i \geq 0 \mid \Ext_A^i(D(A),M) \neq 0 \}$.

\subsection{Spectral sequences}

Spectral sequences are a powerful combinatorial tool for homological algebra. They give a framework to explain long exact sequences involving, for instance, Ext groups and more generally cohomologies and homologies of derived functors. 
For background on spectral sequences, we refer for example to \cite{Rot}. For the convenience of the reader, we recall their definition, elementary properties and some examples to be used throughout this work.

\begin{Def}Let $a\in \mathbb{N}.$
A (cohomology) spectral sequence $\{E_r^{i, j}\colon i, j\in \mathbb{Z}, r\geq a\} $ in an abelian category $\mathcal{A}$ consists of the following data:
\begin{enumerate}[(i)]
\item Pages $E_r$ so that for $r\geq a$, the $r$-page is a collection of objects $\{E_r^{i, j}\colon i, j\in \mathbb{Z}\}$ of $\mathcal{A}$;
\item Boundary maps $d_r^{i, j}\colon E_r^{i, j}\rightarrow E_r^{i+r, j-r+1}$ satisfying the following two conditions \begin{enumerate} 
\item  $d_r^{i, j}\circ d_r^{i-r, j+r-1}=0$;
\item $E_{r+1}^{i, j}=\ker d_r^{i, j}/\im d_r^{i-r, j+r-1}.$
\end{enumerate}
\end{enumerate}
\end{Def}
If for the coordinate $(i, j)$ there exists a degree $r$ such that $E_k^{i, j}=E_r^{i, j}$ for all $k\geq r$, then we write $E_\infty^{i, j}=E_r^{i, j}.$ 
A spectral sequence $\{E_r^{i, j}\colon i, j\in \mathbb{Z}, r\geq a\} $ is called a first-quadrant spectral sequence if all objects $E_r^{i, j}$ are zero unless $i\geq 0$ and $j\geq 0$. 
We say that a first-quadrant spectral sequence converges to $H^*$, written as $E_a^{i, j}\implies H^{i+j}$, if there exists a collection of objects $H^n$ of $\mathcal{A}$ each having a finite filtration
$$0=H_{n+1}^n\subset H_n^n\subset \cdots \subset H_1^n\subset H_0^n=H^n $$ such that $E_\infty^{i, n-i}\cong H_i^n/H_{i+1}^n$ for $0\leq i\leq n$.

Examples of first-quadrant spectral sequences are abundant. For instance, famous examples include the Cartan–Eilenberg Spectral Sequence and the Grothendieck spectral sequence of two left exact functors. For our purposes, the following spectral sequence is crucial.

\begin{Lemma}\label{thespectralsequence}
If $M\in A\m$ and $K$ is a complex with $K^i=0$ for all $i<0$, then there exists a convergent first quadrant spectral sequence
$$ E_2^{i, j}=\Ext_A^i(M, H^j(K))\implies \Ext_A^{i+j}(M, K).$$
\end{Lemma}
\begin{proof}
The result is folklore. It can be found in Lemma 13.21.3 and (15.67.0.1) of the Stacks project (see 
 \citep[\href{https://stacks.math.columbia.edu/tag/0AVG}{Tag 0AVG}]{stacks-project}).
\end{proof}

In many situations, one can directly determine some terms of the limit of the spectral sequence using only information about the second page. A notable example is the cohomology five-term exact sequence (see for instance \citep[Theorem 10.33]{Rot}) which relates $H^1$ and $H^2$ to terms on the second page. This sequence is normally deduced by comparing the filtration on $H^n$ with the differentials on the second page. More generally, when many entries on the second page vanish, a similar strategy often simplifies the computation of the limit terms $H^n$, reducing it to computations involving only terms on the second page.

\begin{Lemma}\label{firstpageequallimit}
	Assume that $E_2^{i, j}\implies H^{i+j}$ is a first quadrant spectral sequence and $E_2^{i, j}=0$ for $i>0$. Then, $E_2^{0, j}\cong H^j$ for every $j\geq 0$.
\end{Lemma}
\begin{proof}
	We claim that $E_s^{i, j}=0$ for $i>0$, $s\geq 2$, $j\geq 0$. We shall proceed by induction on $s$. If $s=2$, the case follows by assumption. Let $s\geq 2$. Then, for any $i>0$, $j\geq 0$, the differentials $d_s^{i, j}$ and $d_s^{i-s, j+s-1}$ must be zero since by induction $E_s^{i, j}=0$. Thus, $E_{s+1}^{i, j}=E_s^{i, j}=0$.

	In particular, $E_\infty^{i, j}=0$ for any $i>0$, $j\geq 0$.
	
	We now claim that $E_s^{0, j}=E_2^{0, j}$ for any $s\geq 2$, $j\geq 0$. We will proceed again by induction. The case $s=2$ is clear. For $s>2$, $j\geq 0$,
	\begin{align}
		E_s^{0, j}=\ker \left( d_{s-1}^{0, j}\colon E_{s-1}^{0, j}\rightarrow E_{s-1}^{s-1, j-s} \right) / \im \left( d_{s-1}^{-s+1, j+s-2}  \right) = \ker \left( d_{s-1}^{0, j}\colon E_{s-1}^{0, j}\rightarrow E_{s-1}^{s-1, j-s} \right),
	\end{align}since $1-s<0$, and thus $d_{s-1}^{-s+1, j+s-2}=0$. By the first claim, $E_{s-1}^{s-1, j-s}=0$. Hence, $E_s^{0, j}=E_{s-1}^{0, j}=E_2^{0, j}$ for any $s\geq 2$, $j\geq 0$. In particular, $E_{\infty}^{0, j}=E_2^{0, j}$ for any $j\geq 0$.
	
	Let $j\geq 0$.	By convergence, we have a filtration for $H^j$ with $E_{\infty}^{0, j}\cong H_0^j/H_1^j=H^j/H_1^j$ for every $j\geq 0$. Furthermore,
	$0=E_\infty^{i, j-i}\cong H_i^j/H_{i+1}^j$, $0<i\leq j$. Thus, $H_i^j=H_{i+1}^j$, $0<i\leq j$. Consequently, $H_1^j=H_{j+1}^j=0$. We conclude $E_2^{0, j}=E_{\infty}^{0, j}\cong H^j$.
\end{proof}

\section{On torsion-free modules and Gorenstein projective modules}
\label{On torsion-free modules and Gorenstein projective modules}

\subsection{Classical dominant dimension} Given $M\in A\m$ we write $\domdim_A M$ to denote the dominant dimension of $M$. Dually, we write $\codomdim_A M:=\domdim_{A^{op}} D(M)$ to denote the codominant dimension of $M$. %The grade of $M$, denoted as $\grade_A M$, is defined as the value $\inf\{i\in \mathbb{N}_0\colon \Ext_A^i(M, A)\neq 0 \}\subset \mathbb{N}_0\cup \{\infty\}$. Dually, one has the cograde of a module, denoted as $\cograde_A M$.
A module $M$ is called \emph{$m$-torsion-free} if $\Ext_A^i(D(A),\tau(M))=0$ for $i=1,...,m,$ where $\tau M$ denotes the Auslander-Reiten translate of $M$. The transpose of $M$ is denoted by $\Tr M$.
We write $\Omega^1 M$ to denote the kernel of $PM\rightarrow M$, where $PM$ is the projective cover of $M$ and $\Omega^{n+1}(M)=\Omega^1(\Omega^n(M))$ for all $n\geq 1$.
A module $N$ is called \emph{$n$-th syzygy} module if $N \cong \Omega^n(M)$
for some module $M$. For $n\geq 1$, we denote by $\Omega^n(\rmod A)$ the subcategory of $\rmod A$ whose modules are isomorphic to direct sums of $n$-th syzygy modules and projective $A$-modules; that is, all modules $M$ that fit into an exact sequence of the form $0\rightarrow M\rightarrow P_{n-1}\rightarrow \cdots \rightarrow P_0$ with each $P_i$ a projective module. By convention, $\Omega^0(\rmod A)=\rmod A$.
$\Dom_m(A)$ denotes the subcategory of modules with dominant dimension at least $m$ and $\TF_m(A)$ the subcategory of modules that are $m$-torsion-free. For algebras of large dominant dimension, those three subcategories all coincide as the next theorem tells us:
\begin{Theorem} \label{connectiontorsionfreedomdim}
Let $A$ be an algebra of dominant dimension at least $n \geq 1$. 
Then we have:
$\Dom_m(A)= \Omega^m(\rmod A)=\TF_m(A)$ for all $1\leq m \leq n$ and $\TF_{n+1}(A)=\Omega^{n+1}(\rmod A)$.
\end{Theorem}
\begin{proof} 
See \cite[Proposition 4]{MV2} and \cite[Proposition 1.6(b), Theorem 0.1]{AR3}.
\end{proof}
The case $n=2$ of the first assertion of Theorem \ref{connectiontorsionfreedomdim} can be traced back to the work of Morita in \cite{zbMATH03271644}.

Recall that the dominant dimension is left-right symmetric, and so the previous result and the coming results are valid for left and also for right modules.

The next theorem is often called Mueller's theorem. It gives a convenient way to calculate the dominant dimension of endomorphism rings of generator-cogenerators.
\begin{Theorem} \label{muellertheo}
Let $A$ be an algebra and $M$ a generator-cogenerator of $\rmod A$. Let $B:=\End_A(M)$. Then the dominant dimension of $B$ is equal to $\inf \{ i \geq 1 | \Ext_A^i(M,M) \neq 0 \} +1$.
\end{Theorem}
\begin{proof}
See \cite[Lemma 3]{zbMATH03248955}.
\end{proof}

In \cite{zbMATH06814513}, the \emph{canonical bimodule} $V:=\Hom_A(D(A),A)$ of a finite-dimensional algebra $A$ was introduced and studied. In particular, the canonical bimodule can be used to characterise properties such as being a Morita or gendo-symmetric algebra, we refer to \cite{zbMATH06814513} for details. This bimodule can also detect whether the underlying algebra has positive dominant dimension.

\begin{Prop}\label{FKYresultforpositivedomdim}
The algebra $A$ has positive dominant dimension if and only if there exists an injective morphism of $A$-bimodules $A\rightarrow \Hom_A(D(A)_A, \Hom_A(V, A_A)_A)$, where $\Hom_A(V, A_A)$ is a right $A$-module through the left $A$-module structure on $V$.
\end{Prop}
\begin{proof}
See \citep[Proposition 4.6]{zbMATH06814513}.
\end{proof}

\subsection{Torsion-free modules and relative dominant dimension}
\label{Torsion-free modules revisited}
%right or left modules?

We shall start by recalling the concept of relative dominant dimension with respect to a module.
Following \cite{cruz2022split}, a module $M$ has \emph{relative dominant dimension} with respect to $Q$, $Q\dd_A M$, at least $n$ if there exists an exact sequence $0\rightarrow M\rightarrow Q_1\rightarrow Q_2\rightarrow \cdots \rightarrow Q_n$ with $Q_i\in \add_A Q$ which remains exact under $\Hom_A(-, Q)$.

This concept (over finite-dimensional algebras) can be characterised by the vanishing of certain $\Ext$-groups.

\begin{Theorem}\label{thm2dot2}
	Let $A$ be a finite-dimensional $K$-algebra and $Q\in A\m$. Let $B$ be the endomorphism algebra $\End_A(Q)^{op}$. For $M\in A\m$, $Q\dd_A M\geq n\geq 2$ if and only if the canonical map $M\rightarrow \Hom_B(\Hom_A(M, Q), Q)$ is an isomorphism and $\Ext_B^i(\Hom_A(M, Q), Q)=0$ for $i=1, \ldots, n-2$.
	
	Moreover, $Q\dd_A M\geq 1$ if and only if the canonical map $M\rightarrow \End_B(\Hom_A(M, Q), Q)$ is injective.
\end{Theorem}
\begin{proof}
	It follows by  \citep[Theorem 3.1.4, Theorem 3.1.1]{cruz2022split}. See also \citep[Lemma 2.3]{cruz2023higher} and Proposition 1.4(1).
\end{proof}
Thus, a module $M$ is torsionless if and only if $A\dd_A M\geq 1$, and it is reflexive if and only if $A\dd_A M\geq 2$.
%Following \cite{AB69}, a module $M$ is called $n$-torsionfree if $\Ext_A^i(\Tr M, A)=0$ for $i=1, \ldots, n$. 
Being $n$-torsion-free for higher $n$ can also be characterised using relative dominant dimension with respect to a module.
 Indeed, the following result can be found in \citep[(2.17)]{AB69} (see also \citep[Theorem 3]{zbMATH01782077}).

 \begin{Theorem}\label{thm2dot3}
 	Let $A$ be a finite-dimensional $K$-algebra and $M\in A\m$. $M$ is an $n$-torsion-free $A$-module if and only if $A\dd_A M\geq n$.
 \end{Theorem}

 We say that a module is \emph{$\infty$-torsion-free} if it is $n$-torsion-free for every $n\in \mathbb{N}$. In view of Theorem \ref{thm2dot3}, these are precisely the modules $M$ satisfying $A\dd_A M=+\infty$.

\subsection{\texorpdfstring{$\infty$} \ -torsion-free modules and Gorenstein projective modules}

%perhaps move theorem 3.1 to beginning of new section
%escrever sobre a seccao que tinha
%escrever section 5 with applications
%ajeitar a notacao, reorganisar e typos para ficar a historia coerente
In this subsection, we summarise some known properties of Gorenstein projective modules using the language of relative dominant dimension and how we can regard them as a special case of $\infty$-torsion-free modules.  

A module $M$ is called \emph{Gorenstein projective} if $\Ext_A^i(M,A)=0=\Ext_A^i(D(A),\tau(M))$ for all $i>0$. Dually $M$ is called \emph{Gorenstein injective} if $D(M)$ is Gorenstein projective.
Recall that a finite-dimensional algebra $A$ is called \emph{Iwanaga-Gorenstein} if the left and right injective dimension of the regular $A$-module coincide and are finite. $\id A_A$ is called the \emph{self-injective dimension} of $A$ or \emph{Gorenstein dimension} when the underlying algebra is Iwanaga-Gorenstein. 
When $A$ is an Iwanaga-Gorenstein algebra, then one also often uses the term maximal Cohen-Macaulay module instead of Gorenstein projective module.

%A module $M$ is called \emph{Gorenstein projective} if $\Ext_A^i(M, A)=\Ext_A^i(\Tr M, A)=0$ for all $i\geq 1$. 
Gorenstein projective modules can be characterised in many different ways (see for example the survey \cite{Chen} and the references therein.):

\begin{Prop}Let $A$ be a finite-dimensional algebra over a field. Let $M\in A\m$. The following assertions are equivalent.\label{prop4dot1}
	\begin{enumerate}[(1)]
		\item $M$ is Gorenstein projective;
		\item $M$ is $n$-torsion-free and $\Tr M$ is $n$-torsion-free for every $n\geq 1$;
		\item $M\in {}^\perp A$ and $M$ is $n$-torsion-free for every $n\geq 1$;
		\item $M\in {}^\perp A$ and $A\dd_A M=+\infty$;
		\item $A\dd_A M=+\infty$ and $A\dd \Tr M=+\infty$;
		\item $A\dd_A M=+\infty$ and $D(A)\cdd_A \tau M=+\infty$.
		\item $M\in {}^\perp A$, $\Hom_A(M, A)\in {}^\perp A_A$ and $M$ is reflexive.
	\end{enumerate}
\end{Prop}
\begin{proof}
	By definition, $M$ is Gorenstein projective if and only if $\Ext_A^i(M, A)=\Ext_A^i(\Tr M, A)=0$ for all $i\geq 1$. The second condition is equivalent to saying that $M$ is $n$-torsion-free for every $n\geq 1$, thus (1) is equivalent to (3). If $M$ is not projective, then $M\cong \Tr \Tr M$ and so the first condition is equivalent to $\Tr M$ being $n$-torsion-free for every $n\geq 1$. Otherwise, $\Tr M$ is zero, which is trivially $n$-torsion-free for every $n\geq 1$. Hence, (1) is equivalent to (2). By Theorem \ref{thm2dot3}, (3) is equivalent to (4), and (2) is equivalent to (5). Since $D(A)\cdd_A D\Tr M=A\dd_A  \Tr M$, (5) is equivalent to (6). By Theorem \ref{thm2dot2}, (4) is equivalent to (7).
\end{proof}

We denote by $A\operatorname{\!- GProj}$ the full subcategory of $A\m$ whose objects consist of Gorenstein projective modules.
Placing extra assumptions on the algebra can simplify the form and the structure of the Gorenstein projective modules. For instance, when $A$ is an Iwanaga-Gorenstein algebra, then every module in ${}^\perp A$ is already Gorenstein-projective. Following \cite{RZ}, finite-dimensional algebras $A$ satisfying the property $A$-$\operatorname{GProj}={}^\perp A$ are called \emph{left weakly Gorenstein} algebras.

Dually, a finite-dimensional algebra $A$ is called \emph{right weakly Gorenstein} if $A^{op}$ is left weakly Gorenstein. The first example of a non right weakly Gorenstein algebra was found in \cite{zbMATH05052963} and non-commutative counterexamples can be found in \cite{RZ}. A systematic construction of non-weakly Gorenstein algebras can be found in \cite{Mar5}.

\section{Bimodule characterisations of the dominant dimension}\label{Bimodule characterisations of the dominant dimension}

In this section, our aim is to study the relation between the relative dominant dimension of the regular module as a bimodule with the absolute dominant dimension, proving the main result.
We begin by recalling the following classical result due to M\"uller:
\begin{Theorem}\label{thm2dot1}
	Let $A$ be a finite-dimensional $K$-algebra. Then
	\begin{align}
		\domdim A^e=\domdim A^{op}=\domdim A=\domdim_{A^e} A=\domdim_{A^e} (A_{A^e}).
	\end{align}
\end{Theorem}
\begin{proof}
	See \citep[Theorems 4 and 8]{zbMATH03248955}.
\end{proof}
The aim is to translate these equalities into properties of $A$ as a module over its enveloping algebra. In particular, we will see that these values characterise torsion-freeness of $A$ as a bimodule and we will get new equivalent characterisations of dominant dimension.

\subsection{Proof of the main result}

 In particular, in this subsection, we want to establish that \begin{equation}
\domdim A=A^e\dd_{A^e} A. \label{bigequality}
\end{equation}

The strategy is as follows:
\begin{itemize}
\item Use the results of \cite{zbMATH03248955} \cite{AR3} and \cite{MV2}  together with the formulation of Theorem \ref{thm2dot3} to deduce $\leq$ in (\ref{bigequality});
\item The proof for the inequality $\geq$ is divided into two cases: 
\begin{itemize}
\item The lower cases $1$ and $2$ are obtained by translating results of \cite{zbMATH06814513} into the language of relative dominant dimension;
\item The higher cases are obtained by employing spectral sequences alongside a recent reformulation of dominant dimension (see \cite{zbMATH07198558}) that can be seen as a continuation of the work \cite{zbMATH06814513}.
\end{itemize}
\end{itemize}

Throughout this subsection, we make the convention that $A$ is a right $A^e$-module and $A^*=\Hom_{A^e}(A, A^e)$ is a left $A^e$-module. Of course, the results in the previous section remain valid for right modules by just considering them over the opposite algebra. From now on, $(-)^*$ denotes the dual $\Hom_{A^e}(-, A^e)$.

\begin{Lemma}\label{maininequalityone}
Let $A$ be a finite-dimensional $K$-algebra. Then $A^e\dd_{A^e} A\geq \domdim A$.
\end{Lemma}
\begin{proof}
 Assume that $\domdim A\geq n$ for some natural number $n$. By Theorem \ref{thm2dot1}, $\domdim A^e\geq n$ and $A_{A^e}\in \Dom_n(A^e)$.  By Theorem \ref{connectiontorsionfreedomdim}, $A$ is $n$-torsion-free over $A^e$. By Theorem \ref{thm2dot3}, $A^e\dd_{A^e} A\geq n$.
\end{proof}
The aim is to show that the inequality in Lemma \ref{maininequalityone} is in fact an equality. The reverse inequality in the lower cases already appears in \cite{zbMATH06814513}, albeit in a different formulation.

\begin{Cor}\label{cor3dot7}
Let $A$ be a finite-dimensional $K$-algebra. Let $i\in \{1, 2\}$. Then, $\domdim A\geq i$ if and only if $A^e\dd_{A^e} A\geq i$ with $A\in \rmod A^e$.
\end{Cor}
\begin{proof}
Suppose that $A^e\dd_{A^e} A\geq i$ for $i\in \{1, 2\}$. By Theorem \ref{thm2dot2}, the map $$A\rightarrow \Hom_{A^e}(\Hom_{A^e}(A, A^e), A^e)$$ is a monomorphism if $i=1$ and an isomorphism if $i=2$. Hence, by Lemma \ref{lem3dot8}, there exists a monomorphism (resp. isomorphism)
\begin{equation}
A\rightarrow \Hom_A(D(A)_A, \Hom_A(\Hom_A(D(A), A), A_A)) \label{iso}
\end{equation} if $i=1$ (resp. $i=2$). By Proposition \ref{FKYresultforpositivedomdim}, $\domdim A\geq 1$. If $i=2$, then it follows by \citep[Theorem 4.7]{zbMATH06814513} and the existence of the isomorphism (\ref{iso}) that $\domdim A\geq 2$.
The converse implication is Lemma \ref{maininequalityone}.
\end{proof}

For the remaining cases, we may assume that the algebra has larger dominant dimension; consequently, it exhibits more favorable properties, and additional methods are available to compute it. We recall the following result.

\begin{Theorem}\label{otherdomdimcharacterisation}
Let $A$ be a finite-dimensional algebra over a field with dominant dimension at least two and let $n$ be a natural number bigger than two. Then, $\domdim A\geq n$ if and only if $$\Ext_{A}^i(D(A)\otimes_A \Hom_A(D(A), A), A_A)=0$$ for $i=1, \ldots, n-2$. 
\end{Theorem}
\begin{proof}
The functors $G:=\Hom_{A^{op}}(\Hom_A(-, A), A)$ and $\Hom_{A^{op}}(D(A)\otimes_A \Hom_A(D(A), A), -)$ are isomorphic by Corollary 3.6 of \citep{zbMATH07198558}. By \citep[Theorem 4.1]{zbMATH06814513}, the dominant dimension of $A$ is at least $n$ if the right derived functor of $G$ vanishes for all degrees, up to $n-2$, that is, $\R^iG=0$ for $i=1, \ldots, n-2$. So the result follows.
\end{proof}

 So to prove the equality $\domdim A=A^e\dd_{A^e} A$ it is enough to show that \begin{equation}A^e\dd_{A^e} (A_{A^e})\geq n \implies \Ext_{A}^i(D(A)\otimes_A \Hom_A(D(A), A), A_A)=0, \ i=1, \ldots, n-2.\label{equalitytoprove} \end{equation}

For cases one and two, the strategy previously used relies on the Tensor-Hom adjunction.
Although, in higher cases, the desired equality may seem like a straightforward application of a higher Tensor-Hom adjunction, this is far from being the case. In fact, the Tensor-Hom adjunction does not generally extend to an Ext-Tor adjunction without further assumptions. Moreover, applying the functor $D(A)\otimes_A -$ to a projective resolution of $\Hom_A(D(A), A)$ as a bimodule does not necessarily yield a 
projective resolution of $D(A)\otimes_A \Hom_A(D(A), A)$.

Indeed as the following example illustrates the vanishing of $\Tor_i^{A^e}(D(A^e), A^*)$ does not imply the vanishing of $\Tor_i^A(D(A), A^*).$

\begin{Example}
Let $Q$ be the quiver $\begin{tikzcd}
& 1 \arrow[dl, "a", swap]& \\
3 \arrow[rr, "b"] & & 2 \arrow[ul, "c", swap] 
\end{tikzcd}$ and $A=KQ/I$ with $I$ generated by $ac$ and $cb$. The algebra $A$ has dominant dimension $3$, so $\Ext_{A^e}^1(\Hom_{A^e}(A, A^e), A^{e})=0$ and $\Ext_A^1(D(A)\otimes_A \Hom_A(D(A), A), A)=0.$ The module $\Hom_A(D(A), A)$ as left $A$-module is isomorphic to $P(1)\oplus P(2)\oplus \top P(1)$, where $P(i)$ is the projective cover of the simple module indexed by $i$. Since $\Ext_A^1(\top P(1), P(3))\neq 0$ it follows that $\Ext_A^1(\Hom_A(D(A), A), A)\neq 0.$ Thus, $\Tor_1^A(D(A), \Hom_A(D(A), A))\neq 0$.
\end{Example}

Spectral sequences remain a very important and useful tool to compute cohomologies and Ext groups. We will use them now to prove (\ref{equalitytoprove}).
By Lemma \ref{lem3dot8},
 \begin{align}
\Ext_{A^e}^i(A^*, A^e)&\cong \Ext_A^i(D(A)\otimesL_A A^*, A_A)\cong \Ext_A^i(D(A_A), D(D(A)\otimesL_A A^*))\\&\cong \Ext_A^i(D(A), \RHom_A(A^*, DDA )\cong \Ext_A^i({}_AD(A), \RHom_A(A^*, A )). \label{eq10}
 \end{align}
 Analogously, 
 \begin{align}
\Ext_A^i(D(A)\otimes_A A^*, A_A)&\cong \Ext_A^i(D(A), D(D(A)\otimes_A A^*))\cong \Ext_A^i(D(A), \Hom_A(A^*, DDA))\\&\cong \Ext_A^i({}_AD(A), \Hom_A(A^*, A)). \label{eq13}
 \end{align}

As injective resolutions of modules are by definition concentrated in non-negative degrees, it follows that $\RHom_A(A^*, A)$ is a complex concentrated in non-negative degrees.
Applying Lemma \ref{thespectralsequence} with $K=\RHom_A(A^*, A)$ and $M={}_AD(A)$ yields the spectral sequence 
\begin{align}
E_2^{i, j}=\Ext_A^i(D(A), H^j(\RHom_A(A^*, A))\implies \Ext_A^{i+j}(D(A), \RHom_A(A^*, A))\label{eq14}
\end{align} Observe that $E_2^{i, j}=\Ext_A^i(D(A), \Ext_A^j(A^*, A))$ and therefore by (\ref{eq13}) the vanishing of $E_2^{i, 0}$ controls the dominant dimension of $A$. We  write $H^i$ to denote $\Ext_A^{i}(D(A), \RHom_A(A^*, A))$.

 \begin{Example}
Assume that $\domdim A\geq 2$. By the cohomology five-term exact sequence (see for instance \citep[Theorem 10.33]{Rot}) there exists an exact sequence
 \begin{align}
0\rightarrow E_2^{1, 0}\rightarrow H^1\rightarrow E_2^{0, 1}\rightarrow E_2^{2, 0}\rightarrow H^2 \label{eq15}
 \end{align}

If $A^e\dd_{A^e} A\geq 3$, then $\Ext_{A^e}^1(A^*, A^e)=0$. By (\ref{eq10}), $H^1=0$. Thus, $E_2^{1, 0}=0$. By the discussion above, this means that $\domdim A\geq 3.$
 \end{Example}

As (\ref{eq15}) illustrates, we need to understand the terms $E_2^{0, j}$. To do this, we will make use of yet another spectral sequence.

\begin{Lemma}\label{lemma3dot11}
Assume that $A^e\dd_{A^e} A\geq n+2$ for some natural number $n$. Let $X$ be a projective-injective $A$-module. Then, $\Hom_A(X, \Ext_A^j(A^*, A))=0$ for $j\in \{1, \ldots, n\}$, where $A^*=\Hom_{A^e}(A, A^e)$.
\end{Lemma}
\begin{proof}
Consider the spectral sequence from (\ref{eq14}) now applied to $X$ instead of $D(A)$:
\begin{align}
S_2^{i, j}=\Ext_A^i(X, \Ext_A^j(A^*, A))\implies \Ext_A^{i+j}(X, \RHom_A(A^*, A)).
\end{align} Observe that $S_2^{i, j}=0$ for all $i>0$ because $X$ is a projective $A$-module. By Lemma \ref{firstpageequallimit}, $S_2^{0, j}\cong \Ext_A^{j}(X, \RHom_A(A^*, A))$ for all $j\geq 0$. Since $A^e\dd_{A^e} A\geq n+2$ we obtain that $\Ext_{A^e}(A^*, A^e)=0$ for $i=1, \ldots, n$. By (\ref{eq10}), $\Ext_A^i(D(A), \RHom_A(A^*, A))=0$ for $i=1, \ldots, n$. Since $X$ is injective it follows that $\Ext_A^i(X, \RHom_A(A^*, A))=0$ for $i=1, \ldots, n$. Therefore, $S_2^{0, j}=0$ for $j=1, \ldots, n$.
\end{proof}

The following result illustrates that the dominant dimension is a lower bound for the cograde of $\Ext_A^j(A^*, A).$
\begin{Cor}\label{lastpiece}
If $\domdim A\geq n+2$ for some natural number $n$, then $$\Ext_A^i(D(A), \Ext_A^j(A^*, A))=0$$ for $i=0, \ldots, n+1$ and $j=1, \ldots, n$, where $A^*$ denotes the module $\Hom_{A^e}(A, A^e)$.
\end{Cor}
\begin{proof}
Let \begin{align}
\cdots \rightarrow P_1\rightarrow P_0\rightarrow D(A)\rightarrow 0 
\end{align} be the minimal projective resolution of $D(A)$. Since $\domdim A\geq n+2$, $D(A)$ has codominant dimension at least $n$, that is, the modules $P_0, \ldots, P_{n+1}$ are projective and injective. Let $j\in \{1, \ldots, n\}$. Fix $M_j:=\Ext_A^j(A^*, A)$.  By Lemma \ref{lemma3dot11}, $\Hom_A(P_i, M_j)=0$ for $i=0, \ldots, n+1$. 
Hence, $\Ext_A^i(D(A), M_j)=0$ for $i=1, \ldots, n+1$ since $\Ext_A^i(D(A), M_j)$ is a quotient of a submodule of $\Hom_A(P_i, M_j)=0$.
\end{proof}

We finally have all the tools to further extend Theorem \ref{thm2dot1}.

\begin{Theorem} \label{mainresult}
Let $A$ be a finite-dimensional $K$-algebra. Then
	\begin{align}
		\domdim A=A^e\dd_{A^e} (A_{A^e})=A^e\dd_{A^e} A.
	\end{align}
\end{Theorem}
\begin{proof}
 First, consider the case where $A$ is a right $A^e$-module. The inequality $\domdim A\leq A^e\dd_{A^e} A$ follows from Lemma  \ref{maininequalityone}. It remains to prove that if $A^e\dd_{A^e} A\geq n$, then $\domdim A\geq n$ for all natural numbers $n$.
As we have seen, if $n$ is one or two, then this statement follows from Corollary \ref{cor3dot7}. We shall prove by induction that $A^e\dd_{A^e} A\geq m+2$ implies $\domdim A\geq m+2$ for every natural number $m$. 
If $m=1$, then the claim follows from the cohomology five-term exact sequence applied to the spectral sequence given in (\ref{eq14}), (see (\ref{eq15})). Consider the first-quadrant spectral sequence $E_2^{i, j}\implies H^{i+j}$ established in (\ref{eq14}).

Assume that $A^e\dd_{A^e} A\geq (m+1)+2=m+3$. By  Theorem \ref{thm2dot2}, and Equation (\ref{eq10}), $H^j=0$ for $j=1, \ldots, m+1$. By induction, we can suppose that $\domdim A\geq m+2$. In particular, the dominant dimension of $A$ is at least two. Recall that $A^*\cong \Hom_A(D(A), A)$. Hence, by Theorem \ref{otherdomdimcharacterisation} and (\ref{eq13}), we obtain that $E_2^{i, 0}=0$ for $i=1, \ldots, m$. 

By Corollary \ref{lastpiece}, $E_2^{i, j}=0$ for $i=0, \ldots, m+1$ and $j=1, \ldots, m$. In particular, $E_2^{i, j}=0$ whenever $i+j=m$ because the spectral sequence is a first-quadrant spectral sequence. Since the element $E_r^{i, j}$ is a quotient of a submodule of $E_{r-1}^{i, j}$ it follows that all elements $E_r^{i, j}=0$ whenever $i+j=m$ for all $r\geq 2$. By definition, $$E_{r+1}^{m+1, 0}=\ker (E_{r}^{m+1, 0}\rightarrow E_{r}^{m+1+r, 1-r})/\im (E_r^{m+1-r, r-1}\rightarrow E_r^{m+1, 0}).$$

Since $1-r<0$ for $r\geq 2$ and $m+1-r+r-1=m$ both terms $E_{r}^{m+1+r, 1-r}$ and $E_r^{m+1-r, r-1}$ are zero. Thus, $E_{r+1}^{m+1, 0}=E_r^{m+1, 0}$ for every $r\geq 2$ and therefore $E_\infty^{m+1, 0}=E_2^{m+1, 0}$.

By convergence, $E_\infty^{m+1, 0}\cong H_{m+1}^{m+1}/H_{m+2}^{m+1}=0$ since $H_{m+1}^{m+1}\subset H^{m+1}=0$. Hence, $E_2^{m+1, 0}=0$ and by Theorem \ref{otherdomdimcharacterisation} and (\ref{eq13}) it follows that $\domdim A\geq m+3.$ The last equality follows then by the left-right symmetry of dominant dimension.
\end{proof}

Observe that because of Theorem \ref{mainresult} we no longer need to be precise whether $A$ is a right or left $A^e$-module.

\subsection{Second formula for dominant dimension using syzygies}

The formula established in Theorem \ref{mainresult} not only provides a new symmetry for dominant dimension, but also provides new structural information on the bimodule $A$.  In fact, this gives us a framework to decide whether $A\in \Omega^n(A^e\m)$ is true using only dominant dimension (see Corollary \ref{cornsyzygydominantdimension} below). The interesting case is when $A$ can be viewed as an $n$-syzygy without being projective as a bimodule.
%Observe that from this new formula, it is also immediate that if $A$ is projective as a bimodule, then $A$ has infinite dominant dimension. 

\begin{Remark}
Actually, $A$ is projective as a bimodule if and only if $A^e$ is semi-simple if and only if $A$ is a separable algebra. The latter is a proper class of semi-simple algebras (see for instance \citep[Proposition 11.8, Theorem 11.11]{SkoYam}).
\end{Remark}

Recall that $J_n(M):=\Tr(\Omega^n(M))$ is the \emph{higher Auslander-Bridger transpose} for $n\in \mathbb{N}_0$.
We denote by $V:=\Hom_A(D(A),A)=\Hom_{A^e}(A,A^e)=A^{*}$ the canonical module, first defined in \cite{zbMATH06814513}.

We now collect our results in the following theorem that gives a new characterisation of the dominant dimension of finite-dimensional algebras:

\begin{Theorem}\label{thm2dot7}
	Let $A$ be a non-separable finite-dimensional algebra over a field and $n\in \mathbb{N}$. Then the following are equivalent:
	\begin{enumerate}[(1)]
		\item $A$ has dominant dimension at least $n$.
		\item $A$ is $n$-torsion-free as a bimodule.
		\item $A^e\dd_{A^e} A\geq n$.
		\item $A\in \Omega^n(A^e\m)$.
	\end{enumerate}
    If $n\geq 2$, then the previous assertions are also equivalent to
    \begin{enumerate}[(i)]
    \item $A\cong \Omega^n(J_{n-2}(V))$ as bimodules.
    \end{enumerate}
\end{Theorem}
\begin{proof}
	By Theorem \ref{thm2dot3}, (2) and (3) are equivalent. By Theorem \ref{mainresult}, (1) is equivalent to (3).  Assume that (4) holds. We claim by induction on $n$ that $A$ is $n$-torsion-free as a bimodule. 
    
    Assume that $n=1$. Then, there exists a monomorphism $A\rightarrow P$ for some $P\in A^e\proj$. Since the injective map $A\rightarrow P\rightarrow P^{**}$ factors through the canonical map $A\rightarrow A^{**}$, this map must be injective. By Theorem \ref{thm2dot2}, $A$ is $1$-torsion-free as a bimodule.

    Assume now that $n>1$.
	Then $A\in \Omega^n(A^e\m)\subset \Omega^{n-1}(A^e\m)$. By induction, $A$ is $(n-1)$-torsion-free as a bimodule. By Theorems \ref{thm2dot3} and \ref{mainresult}, $A^e$ has dominant dimension at least $n-1$. By Theorem \ref{connectiontorsionfreedomdim}, $A$ is $n$-torsion-free as $A^e$-module and so (2) follows. 
Assume that (3) holds. Then by definition $A$ fits into an exact sequence $0\rightarrow A\rightarrow P_1\rightarrow \cdots \rightarrow P_n$ where $P_i$ are projective $A^e$-modules, so (4) is clear.

Now assume that $n\geq 2$. It is clear that (i) implies (4). 
	Assume that (3) holds. By Theorem \ref{thm2dot2}, $V\cong A^*$ and $\Ext_{A^e}^i(\Hom_{A^e}(A, A^e), A^e)=0$ for $i=1, \ldots, n-2$. Let $P_{n-1}\rightarrow P_{n-2}\rightarrow \cdots P_0\rightarrow V\rightarrow 0$ be a minimal projective resolution of $V\cong \Hom_{A^e}(A, A^e)$. Applying the functor $\Hom_{A^e}(-, A^e)$ then yields the exact sequence
	\begin{align}
		0\rightarrow V^*\rightarrow P_0^*\rightarrow \cdots \rightarrow P_{n-1}^*. \label{eq12}
	\end{align} The minimal projective presentation $P_{n-1}\rightarrow P_{n-2}\rightarrow \Omega^{n-2}(V)$, by definition, gives that the cokernel of the map $P_{n-2}^*\rightarrow P_{n-1}^*$ is equal to $\Tr(\Omega^{n-2}(V))$. Since (\ref{eq12}) can be extended to a minimal projective resolution of $\Tr(\Omega^{n-2}(V))$ we obtain that $A\cong A^{**}\cong V^*\cong \Omega^n(\Tr(\Omega^{n-2}(V)))$ as $A^e$-modules.
\end{proof}

\begin{Cor}\label{cornsyzygydominantdimension}
Let $A$ be a finite-dimensional $K$-algebra. Then $$\domdim A=\sup\{n\in \mathbb{N}_0\colon A\in \Omega^n(A^e\m)\}.$$
\end{Cor}
\begin{proof}
If $A$ is separable both sides are $+\infty$. If $A$ is not separable, then the equality follows from Theorem \ref{thm2dot7}.
\end{proof}

%write as Corollary the last equality

\subsection{Revisiting Mueller's theorem and gendo-symmetric algebras} \label{Revisiting Mueller's theorem and gendo-symmetric algebras}

Over the past 15 years, Mueller's characterisation of dominant dimension has gained increasing prominence for both theoretical and practical purposes. In particular, several specialised versions of the characterisation have been developed for specific classes of algebras, for instance for gendo-symmetric algebras \citep[Proposition 3.3]{FanKoe2} and certain classes of Morita algebras \citep[Proposition 4.12]{zbMATH06814513}.

Under the formula established in Theorem \ref{mainresult} we can rewrite Mueller's result (see \citep[Lemma 3]{zbMATH03248955}) as follows:

\begin{Cor}\label{Muellertheoremusingbimodules}
	Let $A$ be a finite-dimensional algebra over a field with dominant dimension at least two. Then, 
	\begin{align*}\domdim A&=\inf\{i\geq 1| \Ext_{A^e}^i(A^*, A^e)\neq 0\}+1\\&=\inf\{i\geq 1| \Hom_{\mathcal{D}(A)}(D(A), \RHom_A(A^*, A )[i])\neq 0\}+1.\end{align*}
\end{Cor}
\begin{proof}
	By Theorem \ref{mainresult}, $\domdim A=A^e\dd_{A^e} A$. By Theorem \ref{thm2dot2} with $Q=A^e$ and $M=A_{A^e}$, the first equality holds. The second equality follows then by Equation (\ref{eq10}).
\end{proof}

 Recall that an algebra $A$ is called \emph{gendo-symmetric algebra}, if it is the endomorphism algebra of a generator over a symmetric algebra $B$ (in the sense that $B \cong D(B)$ as $B$-bimodules).  There are many equivalent characterisations of gendo-symmetric algebras, we refer for example to \cite{FanKoe} and \cite{Mar3}.  For instance, in  \citep[Proposition 3.2]{FanKoe2}, Fang and Koenig characterised gendo-symmetric algebras as the finite-dimensional algebras $A$ with the property that $A$ is isomorphic to $A^*=\Hom_{A^e}(A, A^e)\cong \Hom_A(D(A), A)$ as bimodules. Using this characterisation, we obtain Mueller's theorem for gendo-symmetric algebras established in \citep[Proposition 3.3]{FanKoe2} as a corollary.

\begin{Cor}%[c.f. \citep[Proposition 3.3]{FanKoe2}]
	Let $A$ be a gendo-symmetric algebra. Then the dominant dimension of $A$ equals $\inf \{i \geq 1 | \Ext_A^i(D(A),A) \neq 0 \} +1$.
\end{Cor}
\begin{proof}By definition, $\domdim A\geq 2$ and $A^*\cong A$ as bimodules.  By Corollary \ref{Muellertheoremusingbimodules} and Proposition \ref{restrictingfromenvelopingtoregular}(5), it follows that
	$$\domdim A=\inf\{i\geq 1| \Ext_{A^e}^i(A, A^e)\neq 0\}+1=\inf \{i \geq 1 | \Ext_A^i(D(A),A) \neq 0 \} +1.$$
\end{proof} 

\section{Applications to Hochschild (co)homology}\label{Applications to Hochschild (co)homology}

We apply Theorem \ref{thm2dot7} to give new formulas for the Hochschild cohomology and homology of finite-dimensional algebras with dominant dimension at least two. For the definition and basic properties of the Hochschild homology and cohomology we refer for example to \cite{W}.
Recall that the functor $\tau_{n-1}:=D \Tr \Omega^{n-2}= \tau \Omega^{n-2}$ is called \emph{higher Auslander-Reiten translate}, following \cite{Iya}.
\begin{Theorem} \label{formulahochschild}
Let $A$ be a finite-dimensional $K$-algebra with dominant dimension $n \geq 2$.
\begin{enumerate}
\item We have for the Hochschild homology and $l \geq 1$:
$$HH_l(A) \cong D \Ext_{A^e}^{l+n}(A, \tau_{n-1}(V)).$$
\item We have for the Hochschild cohomology and $l \geq 1$:
$$HH^l(A) \cong \Ext_{A^e}^{l+n}(D(A), \tau_{n-1}(V)).$$

\end{enumerate}
\end{Theorem}
\begin{proof}
By Theorem \ref{mainresult}, if $A$ has dominant dimension $n$, then $A \cong \Omega^n(\Tr(\Omega^{n-2}(V))$.
\begin{enumerate}
\item 
Now $HH_i(A)=\Tor_i^{A^e}(A,A) \cong D \Ext_{A^e}^i (A,D(A))$, using Proposition \ref{torformula} (1).
Using $A \cong \Omega^n(\Tr(\Omega^{n-2}(V))$, we obtain 
\begin{align}
\Ext_{A^e}^i (A,D(A)) &\cong \Ext_{A^e}^i (\Omega^n(\Tr(\Omega^{n-2}(V)),D(A)) \cong \Ext_{A^e}^{i+n} (\Tr(\Omega^{n-2}(V)),D(A)) \\&\cong \Ext_{A^e}^{i+n}(A, \tau_{n-1}(V)).
\end{align}
\item Here 
\begin{align}
HH^i(A)= \Ext_{A^e}^i(A,A) &\cong \Ext_{A^e}^i (\Omega^n(\Tr(\Omega^{n-2}(V)),A) \cong \Ext_{A^e}^{i+n} (\Tr(\Omega^{n-2}(V)),A) \\&\cong \Ext_{A^e}^{i+n}(D(A), \tau_{n-1}(V)).
\end{align}

\end{enumerate}
\end{proof}

Since the formulas are especially nice for gendo-symmetric algebras, which are exactly those algebras with $V \cong A$ as $A^e$-modules, we state this as a corollary:
\begin{Cor}
Let $A$ be a gendo-symmetric algebra with dominant dimension $n \geq 2$.
\begin{enumerate}
\item We have for the Hochschild homology and $l \geq 1$:
$$HH_l(A) \cong D \Ext_{A^e}^{l+n}(A, \tau_{n-1}(A)).$$
\item We have for the Hochschild cohomology and $l \geq 1$:
$$HH^l(A) \cong \Ext_{A^e}^{l+n}(D(A), \tau_{n-1}(A)).$$

\end{enumerate}
\end{Cor}
\begin{proof}
(1) and (2) follow immediately from the previous theorem.
\end{proof}

In the next corollary, we remark that Theorem \ref{formulahochschild} can be used to give vanishing results for the Hochschild homology and cohomology for algebras with dominant dimension at least two.

\begin{Cor}\label{5dot1}
Let $A$ be a finite-dimensional $K$-algebra with dominant dimension $n \geq 2$.
\begin{enumerate}
\item We have $HH_l(A)=0$ for all $l > \pdim_{A^e}(A)-n$.
\item We have $HH_l(A)=0$ for all $l > \idim_{A^e}(\tau_{n-1}(V))-n$.
\item We have $HH^l(A)=0$ for all $l > \pdim_{A^e}(D(A))-n$.
\item We have $HH^l(A)=0$ for all $l > \idim_{A^e}(\tau_{n-1}(V))-n$.
\end{enumerate}
\end{Cor}
\begin{proof}
(1) and (2) follow immediately from Theorem \ref{formulahochschild}(1), while (3) and (4) follow immediately from Theorem \ref{formulahochschild}(2).
\end{proof}

We remark that part (1) can be used to obtain a quick proof that the Hochschild homology of all higher Auslander algebras over an algebraically closed field vanishes in positive degrees since we have $\pdim_{A^e}(A)=\gldim(A)$ by work of Happel in \cite{Ha} if the field is algebraically closed.
While the projective dimension of the bimodule $A$ is known to be equal to the global dimension of the algebra for algebraically closed fields, it seems that the projective dimension of $D(A)$ as a bimodule is not known and there is no homological description for $\pdim_{A^e}(D(A))$ or equivalently $\idim_{A^e}(A)$ in the literature. This motivates us to pose the following question:
\begin{Question}
Let $A$ be a finite-dimensional algebra. Is there a nice homological description of the projective dimension of $D(A)$ as a bimodule, or equivalently of the injective dimension of $A$ as a bimodule?
\end{Question}

We have the following result:
\begin{Prop}\label{prop5dot13}
Let $A$ be a finite-dimensional algebra with enveloping algebra $A^e$.
Then 
$$ \idim_{A^e}(A) \geq \sup \{ \idim_A (A_A), \idim_A ({}_A A) \}.$$
\end{Prop}
\begin{proof}
 If $\idim_{A^e}(A)$ is infinite, then there is nothing to prove. Observe that every injective $A^e$-module is also an injective left $A$-module and an injective right $A$-module. So, the minimal injective resolution of $A$ as $A^e$-module is by restriction also an injective resolution of $A$ as a  left $A$-module (resp. as a right $A$-module). 
\end{proof}
\begin{Cor}Let $A$ be a finite-dimensional algebra over an algebraically closed field.
Then, $A$ has finite global dimension if and only if  $\idim_{A^e}(A)$ is finite. \label{finitenessglobalusingregularasbimodule}
\end{Cor}
\begin{proof}
	 Assume that $A$ has finite global dimension. By the last equation in \cite{zbMATH03114458}, $A^e$ has finite global dimension equal to $2 \gldim A$ and so $\idim_{A^e}(A)$ is finite. Conversely, assume that $\idim_{A^e}(A)$ is finite. By Proposition \ref{prop5dot13}, $A$ is Iwanaga-Gorenstein. Thus, $A^e$ is also Iwanaga-Gorenstein (see for instance \cite[Proposition 2.2 (a)]{ARCohen}). %over algebraically closed fields one can use again zbMATH03114458
	Thus, $\idim_{A^e}(A)$ being finite implies that $\pdim_{A^e} (A)$ is finite as over an Iwanaga-Gorenstein algebra $B$ with $B$-module $N$ one has $\pdim N< \infty$ if and only if $\idim N < \infty$. Now $\pdim_{A^e} (A)$ is precisely the global dimension of $A$ by the lemma in 1.5 of \cite{Ha}, so the result follows.

\end{proof}

We observe now that the injective dimension of $A$ is always equal to the projective dimension of $D(A)$ as $A$-bimodules, or stated equivalently:
\begin{Prop}\label{injdimasbimoduleleftrightsymm}
Let $A$ be a finite-dimensional algebra over an algebraically closed field. Then, 
$$\injdim_{A^e} ({}_{A^e}A)= \injdim_{A^e} (A_{A^e}).$$
\end{Prop}
\begin{proof}
We saw in Corollary \ref{finitenessglobalusingregularasbimodule}, that both terms are infinite if and only if $A$ has infinite global dimension. So, we can assume without loss of generality that $A$ has finite global dimension. Hence,
\begin{align*}
\injdim_{A^e} ({}_{A^e}A)&=\sup \{n\geq 0\colon \Ext_{A^e}^n(D(A^e), ({}_{A^e}A))\neq 0 \}=\sup \{n\geq 0\colon \Ext_{A^e}^n(D(A), A^e_{A^e})\neq 0 \} \\
\injdim_{A^e} (A_{A^e})&= \sup \{n\geq 0\colon \Ext_{A^e}^n(D(A^e), A_{A^e})\neq 0 \}=\sup \{n\geq 0\colon \Ext_{A^e}^n(D(A), {}_{A^e}A^e)\neq 0 \}.
\end{align*} 
By Proposition \ref{restrictingfromenvelopingtoregular}, we obtain $\Ext_{A^e}^n(D(A), A^e_{A^e})\cong \Ext_A^n(D(A)\otimesL_A (D(A)), {}_A A)$ and
\begin{align}
\Ext_{A^e}^n(D(A), {}_{A^e} A^e)&\cong \Ext_A^n(D(A)\otimesL_A (D(A))_A, A_A) \\&\cong \Ext_A^n(D(A_A), D(D(A)\otimesL_A (D(A))_A)) \label{eq29again}\\
&\cong \Hom_{\mathcal{D}(A)}(D(A_A)[-n], \RHom_A(D(A), {}_AA)) \label{eq30}\\&\cong \Hom_{\mathcal{D}(A)}(D(A)\otimesL_A D(A_A)[-n], {}_AA)
\end{align}
  Here $\RHom_A(D(A), A)$ is again in $\mathcal{D}(A)$ by the action of $A$ on $D(A)=(DA)_A$ as a right $A$-module. That is, the module in the first entry of the derived Hom complex of Equation (\ref{eq30}) is the module in the second entry of the derived tensor product appearing in Equation (\ref{eq29again}). Thus, in the last line, $D(A)\otimesL_A D(A_A)\in \mathcal{D}(A)$ by the bimodule structure of $D(A)$ in the first entry. So, $\Ext_{A^e}^n(D(A), A^e_{A^e})\cong \Ext_A^n(D(A)\otimesL_A (D(A))_A, {}_A A) \cong \Ext_{A^e}^n(D(A), {}_{A^e} A^e) $ for every $n\in \mathbb{N}$. 
\end{proof}

The argument used in the previous proposition also implies that the grade of $D(A)$ is equal to the cograde of $A$ not only as modules but also as bimodules. 

\begin{Cor}
Let $A$ be a finite-dimensional algebra over a field. Then,
\begin{align*}
\grade_{A^e}(D(A))=\cograde_{A^e}(A).
\end{align*}
\end{Cor}
\begin{proof}
By the proof of Proposition \ref{injdimasbimoduleleftrightsymm}, it follows that $\Ext_{A^e}^n(D(A), {}_{A^e} A^e)\cong \Ext_{A^e}^n(D(A^e), {}_{A^e} A).$ Hence, \begin{align*}\grade_{A^e}(D(A))&=\inf\{n\geq 0\colon \Ext_{A^e}^n(D(A), {}_{A^e} A^e)\neq 0 \}=\inf\{n\geq 0\colon \Ext_{A^e}^n(D(A^e), {}_{A^e} A)\neq 0 \}\\&=\cograde_{A^e}(A). \end{align*}
\end{proof}

%\begin{question}
%Let $A$ be a finite dimensional algebra. Is the injective dimension of $A$ always equal to the projective dimension of $D(A)$ as $A$-bimodules?
%\end{question}

%We saw in the previous corollary that the answer is yes when $A$ has infinite global dimension. 
%Note that when $A$ has finite global dimension, then we have
%$\pdim D(A)=\sup \{ l \geq 0 \mid \Ext_{A^e}^l(D(A),A^e) \neq 0 \}$
%and $\idim A=\sup \{ l \geq 0 \mid \Ext_{A^e}^l(D(A^e),A) \neq 0 \}$.
%So we can as more generally do we have for every $l \geq 0$:
%$\Ext_{A^e}^l(D(A),A^e) \cong \Ext_{A^e}^l(D(A^e),A)?$
%So far QPA found no counterexample, even for algebras of infinite global dimension?! This would imply that the grade of $D(A)$ as a bimodule would be equal to the cograde of $A$ as a bimodule.
%Motivation: We have $grade_{A^e}(A)=cograde_{A^e}(D(A)).$

We especially note the following corollary:
\begin{Cor}\label{boundsforinjdimasbimodule}
Let $A$ have global dimension $g$ and let the field $K$ be algebraically closed.
Then $g \leq \idim_{A^e}(A) \leq 2g$.
\end{Cor}
\begin{proof} We already saw $g \leq \idim_{A^e}(A)$ in Proposition \ref{prop5dot13}.
We have $\gldim A^e=2 \gldim A$ by the last equation in \cite{zbMATH03114458}.
Thus, $\idim_{A^e}(A) \leq 2g$.
\end{proof}

We give two examples that show that the bounds in the previous corollary are optimal.
\begin{example}
The hereditary algebra $KQ$ with $Q$ given by 
% https://q.uiver.app/#q=WzAsMixbMCwwLCIxIl0sWzEsMCwiMiJdLFswLDFdXQ==
\[\begin{tikzcd}
	1 & 2
	\arrow[from=1-1, to=1-2]
\end{tikzcd}\]
has global dimension 1 and $D(A)$ has projective dimension 1 as a bimodule. 
%For the Nakayama algebra with Kupisch series [2,3] the bimodule $D(A)$ has projective dimension 4, while the algebra has global dimension 2.
\end{example}

\begin{example}
	Let $\mathcal{A}_n$ be the bound quiver algebra of a representation-finite block of a Schur algebra with $n$ simples.
	See for example \cite{zbMATH00509864} that $\mathcal{A}_n$ is the algebra $KQ/I$ where $Q$ is the quiver
	\[
	\begin{tikzcd}[every arrow/.append style={bend left}]
		1 \arrow[r, "a_1"]  & 2 \arrow[l, "b_1"] \arrow[r, "a_2"] & 3 \arrow[l, "b_2"] \arrow[r, "a_3"] & 
		4 \arrow[l, "b_3"] & \cdots & m-1 \arrow{r}[above]{a_{m-1}} & m \arrow{l}[below]{b_{m-1}}
	\end{tikzcd}, \quad m\geq 1,
	\]  and the ideal $I$ is generated by the relations $b_{m-1} a_{m-1}, \ b_{i-1}a_{i-1}-a_i b_i,\ a_{i-1}a_i, \ b_i b_{i-1}$ for  $i=2, \ldots, m-1$.
\end{example}By Proposition 4.1 of \cite{zbMATH05124581}, 
\begin{align}
		\dim_k HH^i(\mathcal{A}_n)=\begin{cases}
			n, & i=0 \\
			1, & 1\leq i \leq 2(n-1)\\
			0, & i\geq 2n-1
		\end{cases}.
\end{align}
Observe that $\mathcal{A}_n$ satisfies $\domdim \mathcal{A}_n=2(n-1)=\gldim \mathcal{A}_n$. By Corollary \ref{5dot1}, the Hochschild homology vanishes in every degree and $\pdim_{\mathcal{A}_n^e} D(\mathcal{A}_n)\geq 4(n-1).$ Since $4(n-1)=\gldim \mathcal{A}_n^e$ we get $\pdim_{\mathcal{A}_n^e} D(\mathcal{A}_n)= 4(n-1)= 2\gldim \mathcal{A}_n$.

The next example shows that in general $\idim_{A^e}(A)$ can be strictly larger than $\gldim A$ and strictly smaller than $2 \gldim A$.

\begin{example}
Let $A$ be the $K$-algebra $KQ/I$, where $Q$ is the quiver 
% https://q.uiver.app/#q=WzAsNCxbMCwwLCIxIl0sWzEsMCwiMiJdLFsyLDAsIjMiXSxbMywwLCI0Il0sWzAsMSwiYSJdLFsxLDIsImIiXSxbMiwzLCJjIl1d
\[\begin{tikzcd}
	1 & 2 & 3 & 4
	\arrow["a", from=1-1, to=1-2]
	\arrow["b", from=1-2, to=1-3]
	\arrow["c", from=1-3, to=1-4]
\end{tikzcd}\] and $I$ is the ideal generated by $abc$. Then $A$
has global dimension 2 and $\idim_{A^e}(A)=3.$

\end{example}

 Motivated by some more experiments and examples with gendo-symmetric algebras we pose the following question:
\begin{question}
Let $A$ be a gendo-symmetric algebra of finite global dimension $g$. Do we have $\idim_{A^e}(A)=2g$?
\end{question}

\section{Applications related to homological conjectures}\label{Applications related to homological conjectures}

The following conjecture is known as the Nakayama conjecture and it was stated for the first time by Nakayama in 1958, see \cite{Nak}.
\begin{Conjecture}[Nakayama conjecture] \label{Nakayama conjecture}
A finite-dimensional algebra $A$ is self-injective if and only if $A$ has infinite dominant dimension.
\end{Conjecture} Actually, in the first appearance, the Nakayama conjecture was formulated using bimodules, but due to Mueller's result (Theorem \ref{thm2dot1}), it is not necessary to formulate it using dominant dimension of bimodules.  
In this section, we advocate that the bimodule approach might provide new insights into  the Nakayama conjecture using bimodules and their potential relations with Gorenstein homological algebra.

\subsection{Detecting infinite dominant dimension in the \texorpdfstring{$\mho$-} \ quiver}
Following \cite{RZ}, for a module $M$ one defines $\mho(M)$ as the cokernel of a minimal left $\add(A)$-approximation of $M$.

The \emph{$\mho$-quiver} of a finite-dimensional algebra $A$ is defined as the quiver having vertices the isomorphism classes of indecomposable non-projective $A$-modules $[X]$ and there is an arrow $[\mho(X)] \rightarrow [X]$ for any torsionless indecomposable non-projective module $X$.

We collect some results on $\mho$.
\begin{Prop} \label{ringelzhangpropo}
Let $A$ be a finite-dimensional algebra. Let $M$ be an indecomposable, non-projective $A$-module.
\begin{enumerate}
\item $\mho^k \cong \Tr \Omega^k \Tr$ for $k \geq 1$.
\item If $M$ is torsionless, $\mho(M)$ is indecomposable and not projective and $\Omega( \mho(M)) \cong M$.
\item $[M]$ is the start of a path of length $t \geq 1$ in the $\mho$-quiver if and only if $\Ext_A^i(M,A)=0$ for $i=1,...,t$.
\item $[M]$ is the end of a path of length $t \geq 1$ if and only if $M$ is $t$-torsion-free.

\end{enumerate}

\end{Prop}
\begin{proof}
\begin{enumerate}
\item See the lemma in 4.4. of \cite{RZ}.
\item This is a consequence of Lemma 3.3. in \cite{RZ}.
\item See (1) of the theorem in Section 1.5. of \cite{RZ}.
\item See (2) of the theorem in Section 1.5. of \cite{RZ}.
\end{enumerate}
\end{proof}
Our main result gives a new point of view on the Nakayama conjecture using the $\mho$-quiver of Ringel and Zhang, which we state as a corollary of our results.

\begin{Cor}
The following are equivalent for a finite-dimensional algebra $A$.
\begin{enumerate}
\item $A$ has infinite dominant dimension.
\item $A$ is $\infty$-torsion-free as an $A$-bimodule.
\item There are paths of arbitrarily large length ending in $A$ in the $\mho$-quiver of $A^e$.

\end{enumerate}

\end{Cor} 
\begin{proof}
The equivalence of (1) and (2) follows directly by Theorem \ref{mainresult}, while the equivalence of (2) and (3) is a consequence of Proposition \ref{ringelzhangpropo} (4).
\end{proof}

\subsection{Gorenstein bimodule conjecture}
Self-injective algebras are the Iwanaga-Gorenstein algebras with Gorenstein dimension zero. So an algebra is self-injective if and only if every module is Gorenstein projective. By our main theorem, if $A$ is Gorenstein projective as a bimodule, then the dominant dimension of $A$ is infinite.

Following the by now well-known meta-theorem: \newline
	\textit{Every result in homological algebra has a counterpart in Gorenstein homological algebra} \newline
 we should expect a Gorenstein counterpart to the equality $\pdim_{A^e} A=\gldim A$ for finite-dimensional algebras $A$. Within this expectation, we give the following new Gorenstein homological conjecture, called the Gorenstein bimodule conjecture. Further motivation for this conjecture will be given afterwards by relating it to the Nakayama and Tachikawa conjectures.
\begin{Conjecture}[Gorenstein bimodule conjecture] \label{gorbimodconjecture}
	A finite-dimensional algebra $A$ is self-injective if and only if $A$ as a bimodule is Gorenstein projective.
\end{Conjecture}
The conjecture is known to be true in case $A$ is Iwanaga-Gorenstein by results of Shen in \cite{S} and more generally when $A$ is left weakly Gorenstein by results in \cite{zbMATH07040786}.
Indeed, in both situations if $A\in A^e\operatorname{\!- GProj}$, then $D(A)\in {}^\perp {}_AA$ and thus $D(A)$ is a Gorenstein projective module over $A$. Hence $D(A)$ is then torsionless and thus a submodule of a projective module forcing $A$ to be self-injective.  Below in Proposition \ref{prop6dot8} we give an easy proof that the Gorenstein bimodule conjecture is true for right weakly Gorenstein algebras.

In \cite{zbMATH07971443}, the authors introduced and characterised algebras with $n$-torsion-free Auslander-Reiten sequences. These are precisely the algebras affording the property: $\domdim A\geq n$ and $\domdim \End_A(A\oplus D(A))\geq n+2$. There, the authors conjectured that a finite-dimensional algebra has $n$-torsion-free Auslander-Reiten sequences for every $n\in \mathbb{N}$ if and only if it is self-injective (we refer to \cite{zbMATH07971443} for more details and definitions). The following theorem shows that Conjecture \ref{gorbimodconjecture} is equivalent to Conjecture 3.9 of \cite{zbMATH07971443}.

\begin{Theorem}\label{thm6dot4}
	Let $A$ be a finite-dimensional algebra over a field. The following assertions are equivalent.
	\begin{enumerate}[(1)]
		\item $A$ is Gorenstein projective as a bimodule;
		\item $\domdim A=+\infty$ and $D(A)\in {}^\perp {}_AA$;
		\item $\domdim A=+\infty$ and $D(A)\in {}^\perp A_A$;
		\item $\domdim_{A^e} A=\domdim_{A^e} \Hom_{A^e}(A, A^e)=+\infty$;
		\item $\domdim_{A^e} A=\codomdim_{A^e} \tau A=+\infty$;
		\item $\domdim A=\domdim \End_A(A\oplus D(A))=+\infty$;
		\item $A$ has $n$-torsion-free Auslander-Reiten sequences for every $n\geq 1$;
		\item $A\in {}^\perp A^e$ and $A^e\dd_{A^e} A=+\infty$.
	\end{enumerate}
\end{Theorem}
\begin{proof}
	The equivalence between (6) and (7) follows from \citep[Theorem 3.2]{zbMATH07971443}.
	
	From $\Ext_A^i(D(A_A), {}_A A)\cong \Ext_A^i(D({}_A A), DDA_A)\cong \Ext_A^i(D({}_A A), A_A)$ follows the equivalence between (2) and (3). By Proposition \ref{prop4dot1}, (1) is equivalent to (8). By Proposition \ref{restrictingfromenvelopingtoregular}(5), the statements (8) and (3) are equivalent. By Theorem \ref{muellertheo}, $D(A)\in {}^\perp A$ if and only if $\domdim \End_A(A\oplus D(A))=+\infty$. By Theorem \ref{thm2dot1}, $\domdim_{A^e} A=\domdim A=\domdim A^e$. So the statements (2) and (3) are equivalent to (6). By Theorem \ref{thm2dot1}, we also obtain that $\domdim_{A^e} A^*=+\infty$ if and only if $A^e\dd_{A^e} A^*=+\infty$ under the assumption $\domdim A=+\infty$. By Theorem \ref{thm2dot2}, $A^e\dd_{A^e} A^*=+\infty$ if and only if $\Ext_{A^e}^i(A, A^e)\cong \Ext_{A^e}^i(A^{**}, A^e)=0$ for all $i\geq 1$ under the assumption $\domdim A=+\infty$. So (4) is equivalent to (8). Assuming that $\domdim A^e=+\infty$ we obtain $\domdim_{A^e} A^* =\domdim_{A^e} \Tr A+ 2$ using the exact sequence $0\rightarrow A^*\rightarrow P_0^*\rightarrow P_1^*\rightarrow \Tr A\rightarrow 0$. It follows that (4) and (5) are equivalent.
\end{proof}

The next conjecture is called the first Tachikawa conjecture and was first stated in the book \cite{Tachikawabook}. 
\begin{Conjecture} \label{Tachikawa conjecture}
	A finite-dimensional algebra $A$ (over a field) is self-injective if and only if $\Ext_A^i(D(A),A)=0$ for all $i \geq 1$.
\end{Conjecture}

It is known that the truth of the Nakayama conjecture for all algebras $A$ would imply the Tachikawa conjecture, but it is not known whether the truth of the Nakayama conjecture for a fixed algebra implies also that the first Tachikawa conjecture is true for this fixed algebra.
Recall that the \emph{finitistic dimension} of a finite-dimensional algebra $A$ is defined as the supremum of all projective dimensions of modules having finite projective dimension. The finitistic dimension conjecture predicts that all finite-dimensional algebras have finite finitistic dimension. It is known that the finitistic dimension conjecture implies the Nakayama conjecture for a given algebra $A$. See, for example, \cite{Yam} for further discussion of these conjectures and their connections to other homological conjectures.

As an application of our main result of this article, we can relate the Gorenstein bimodule conjecture to the Nakayama and first Tachikawa conjectures:
\begin{Theorem} \label{gorbimodtheorem}
	Let $A$ be a finite-dimensional algebra.
	Then the following are equivalent for $A$:
	\begin{enumerate}
		\item The Gorenstein bimodule conjecture holds for $A$.
		\item The Nakayama conjecture or the first Tachikawa conjecture holds for $A$.
		
	\end{enumerate}
    \end{Theorem}
\begin{proof}
 It follows directly by (2) of Theorem \ref{thm6dot4}.
\end{proof}

As suggested by the referee, right weakly Gorenstein algebras are precisely those algebras for which the left condition in (4) of Proposition \ref{prop4dot1} can be omitted. In particular, the Nakayama conjecture and the Gorenstein bimodule conjecture are true for right weakly Gorenstein algebras.

%escrever a nova prova combinando 6.8 e 6.9 depois escrever comentarios

\begin{Prop}\label{prop6dot8}
Let $A$ be a finite-dimensional algebra. Then, $A$ is right weakly Gorenstein if and only if \begin{align}A\operatorname{\!- GProj}=\{X\in A\m\colon A\dd_A X=+\infty\}. \label{eq35}\end{align}
Moreover, a right weakly Gorenstein algebra $A$ is self-injective if and only if $\domdim A=+\infty$.
\end{Prop}
\begin{proof}
Assume that $A$ is right weakly Gorenstein. By Proposition \ref{prop4dot1}, every Gorenstein projective module $M$ satisfies $A\dd_A M=+\infty$. Let $M$ be a module satisfying $A\dd_A M=+\infty$. By Theorem \ref{thm2dot2}, $M^*\in {}^\perp A_A$ and $M\cong M^{**}$, where $(-)^*=\Hom_A(-, A)$. Since $A^{op}$ is left weakly Gorenstein, $M^*$ is Gorenstein projective and so $M\cong M^{**}\in {}^\perp A$. Thus, $M\in A\operatorname{\!- GProj}$.

Conversely assume that (\ref{eq35}) holds. Let $Y\in {}^\perp A_A$. Thus, $\Ext_A^i(\Omega^2 Y, A_A)\cong \Ext_A^{i+2}(Y, A_A)=0$ for $i>0$. Let $P_1\rightarrow P_0\rightarrow Y$ be the minimal projective presentation of $Y$. Then, the exact sequence \begin{align}
0\rightarrow \Omega^2(Y)\rightarrow P_1\rightarrow P_0\rightarrow Y \rightarrow 0\label{eq32}
\end{align} remains exact under $\Hom_A(-, A)$ and so $A_A\dd_{A^{op}} \Omega^2(Y)\geq 2.$ By Theorem \ref{thm2dot3}, $X:=\Omega^2(Y)$ is reflexive. So $X^*$ is reflexive and $(X^*)^*\cong X\in {}^\perp A_A$. Thus ${}_AA\dd_A X^*=+\infty$. By (\ref{eq35}), $X^*\in A\operatorname{\!- GProj}$ and thus $X^*\in {}^\perp {}_AA$. Thus, $\Omega^2(Y)=X\in A^{op}\operatorname{\!- GProj}$. By (\ref{eq32}) and Proposition \ref{prop4dot1}, $2+A_A\dd_{A^{op}} Y=A_A\dd_{A^{op}} \Omega^2(Y)=+\infty$.
Again by Proposition \ref{prop4dot1}(4), it follows that $Y\in A^{op}\operatorname{\!- GProj}$. That is, $A^{op}$ is left weakly Gorenstein, so the first statement holds.

To show the second statement, assume that $A$ is right weakly Gorenstein and suppose that $\domdim A=+\infty$. Consider the injective copresentation of $A$:
\begin{align}
	0\rightarrow A\rightarrow I_0(A)\rightarrow \Omega^{-1}(A)\rightarrow 0. \label{eq34}
\end{align} Since $I_0(A)$ is projective-injective and $\domdim A=+\infty$, it follows that $\domdim_A  \Omega^{-1}(A)=+\infty$. By Theorems \ref{connectiontorsionfreedomdim} and \ref{thm2dot3}, $A\dd_A \Omega^{-1}(A)=+\infty$. By assumption, $\Omega^{-1}(A)$ is Gorenstein-projective and so $\Omega^{-1}(A)\in {}^\perp A$. Thus, the exact sequence (\ref{eq34}) splits and so $A$ is self-injective.
\end{proof}

Since dominant dimension is left-right symmetric, it follows from Proposition \ref{prop6dot8} that the Nakayama conjecture is also true for left weakly Gorenstein algebras.

\begin{Remark}
If $\injdim_A {}_A A\leq n<+\infty$, then $A$ is right weakly Gorenstein. This is known by \citep[4.8]{RZ} and \cite{zbMATH02174886}. Using  Proposition \ref{prop6dot8} the argument is as follows: Let $X\in A\m$ with $A\dd_A X=+\infty$, then applying $\Hom_A(-, A)$ to an exact sequence giving $A\dd_A X=+\infty$ we get by dimension shift that $X\in {}^\perp A$ because $A$ has finite injective dimension. So \eqref{eq35} holds. In particular, $A\operatorname{\!- GProj}\subset \Omega^n(A\m)$.
 Moreover, if $A$ is, in addition, left weakly Gorenstein, then since $\Omega^n(A\m)\subset {}^\perp A$ we obtain that $\Omega^n(A\m)=A\operatorname{\!- GProj}$ which implies that $A$ is Iwanaga-Gorenstein (see for instance \cite{Chen}).
	\end{Remark}
 
\begin{Cor}
	Let $A$ be an algebra of finite finitistic dimension, then $A$ satisfies the Gorenstein bimodule conjecture.
	
\end{Cor}
\begin{proof} 
	Since having finite finitistic dimension implies that $A$ satisfies the Nakayama conjecture, 
	by Theorem \ref{gorbimodtheorem} $A$ also satisfies the Gorenstein bimodule conjecture. 
\end{proof}
In particular, the Gorenstein bimodule conjecture is true for Iwanaga-Gorenstein, monomial or local algebras since all such algebras have finite finitistic dimension, and as we have seen it is also true for left weakly Gorenstein algebras and right weakly Gorenstein algebras.

For gendo-symmetric algebras the Gorenstein bimodule conjecture is even equivalent to the Nakayama conjecture, we state this as a corollary in the following form:

\begin{Cor}
	Let $A$ be a gendo-symmetric algebra.
	Then the following are equivalent:
	\begin{enumerate}
		\item $A$ is Gorenstein projective as an $A$-bimodule.
		\item $A$ has infinite dominant dimension.
		
	\end{enumerate}
	\label{cor6::10}
\end{Cor}
\begin{proof}
	 This follows directly from Theorem \ref{thm6dot4} together with the fact that for gendo-symmetric algebras $A$ is isomorphic to  $A^*$ as $A^e$-modules.
\end{proof}

This equivalence exhibited in Corollary \ref{cor6::10} can be extended to a bigger class of algebras that contains gendo-symmetric algebras. Indeed, we can weaken the condition of $A^*$ being isomorphic to the regular module to just requiring that $A^*$ is Gorenstein projective as an $A$-module.

\begin{Cor}
Let $A$ be a finite-dimensional algebra over an algebraically closed field such that $A^*=\Hom_{A^e}(A, A^e)$ is Gorenstein projective as a right and as a left $A$-module.
Then $A$ has infinite dominant dimension if and only if $A$ is Gorenstein projective as a bimodule.
\end{Cor}
\begin{proof}
By Theorem \ref{thm6dot4}, $A$ has infinite dominant dimension whenever $A$ is Gorenstein projective as a bimodule.

Conversely, assume that $\domdim A=+\infty$. By Proposition \ref{prop4dot1}, $A\dd_A ({}_AA^*)=+\infty$  and $A\dd_A (A^*_{A})=+\infty.$ Since dominant dimension is left-right symmetric, Theorem \ref{connectiontorsionfreedomdim} yields that $\domdim_A ({}_AA^*)=\domdim_A (A^*_A)=+\infty$.
By Lemma 7 of \cite{zbMATH03248955}, it follows that $\Ext_A^i(S, {}_AA^*)=\Ext_{A}^i(S', A^*_A)=0$ for every $i\geq 0$ and $S$ simple left $A$-module not in the socle of ${}_AA$ and $S'$ simple right $A$-module not in the socle of $A_A$.
Let $W$ be a simple $A^e$-module which is not in the socle of $A^e$. Since the underlying field is algebraically closed, then $W\cong S_1\otimes_K S_2$ for $S_1$ a simple left $A$-module and $S_2$ a simple right $A$-module. Moreover, either $S_1$ is not in the socle of ${}_A A$ or $S_2$ is not in the socle of $A_A$. Assume that the latter occurs. 

By \citep[Theorem 2.8a, p.167]{MR77480}, $\Ext_{A^e}^i(S_1\otimes_K S_2, A^*)\cong \Ext_{K\otimes_K A}^i(S_1, \Hom_A(S_2, A^*_A))=0$ since $\Ext_A^{n> 0}(S_2, A^*)=0=\Tor_{n>0}^K(S_2, S_1)$.

By Lemma 7 of \cite{zbMATH03248955}, we obtain that $\domdim_{A^e} A^*=+\infty$.
 By Theorem \ref{thm2dot1}, $\domdim_{A^e} A=\domdim A=+\infty$. So, by Theorem \ref{thm6dot4}, the result follows.
\end{proof}

We conclude by noting that a finite-dimensional algebra $A$ is self-injective if and only if the regular module $A$ is Gorenstein injective as a bimodule.

\begin{Remark}
Let $A$ be a finite-dimensional algebra over a field. Then, $A$ is self-injective if and only if $D(A)$ is Gorenstein projective as a bimodule.
\end{Remark}
\begin{proof}
If $A$ is self-injective, then $A^e$ is also self-injective since the tensor product of projective-injective modules is again projective-injective. Since over self-injective algebras every module is Gorenstein-projective, it is clear that $D(A)$ is Gorenstein projective as a bimodule.

Conversely, suppose that $D(A)$ is Gorenstein projective as a bimodule. By Proposition \ref{prop4dot1}, $D(A)$ is, in particular, $1$-torsion-free (that is, torsionless) as a bimodule. 

Recall that, in general, a finitely generated module $M$ is torsionless if and only if $M$ can be embedded into a projective module. Indeed, if $M$ is torsionless, then $M\hookrightarrow M^{**}$. Since $M^*$ is finitely generated, $M^{**}$ can be embedded into a direct sum $(A^s)^*\cong A^s$. The other direction follows from the fact that any injective map $M\rightarrow A^s\rightarrow (A^s)^{**}$  factors through $M\rightarrow M^{**}$. 

So, in our case, we obtain that $D(A)$ can be embedded into a projective module, and so by being injective it must also be projective. Since $A$ is a finite-dimensional algebra it follows that $A$ is self-injective.
\end{proof}

\section{Further questions and a conjecture}\label{Further questions}
Motivated by the results in our article and some small experiments, we give three questions here for general noetherian rings. 
\begin{question}
Let $R$ be a two-sided noetherian ring.
\begin{enumerate}
\item (Gorenstein bimodule conjecture for noetherian rings) Is $R$ as a bimodule Gorenstein projective if and only if $R$ is self-injective?
\item (Nakayama conjecture for noetherian rings) Is $R$ an $n$-syzygy module as a bimodule for every $n \geq 1$ if and only if $R$ is self-injective?
\item What is the injective dimension of the regular module $R$ as a bimodule?
\end{enumerate}
\end{question}
It would be especially interesting to look at those questions when $R$ is additionally commutative. 
We also add the following conjecture:
\begin{conjecture}
Let $A$ and $B$ be two finite dimensional $K$-algebras that are derived equivalent. Assume that $A$ is Gorenstein projective as an $A$-bimodule. Then also $B$ is Gorenstein projective as a $B$-bimodule.
\end{conjecture}
By Theorem \ref{gorbimodtheorem} this conjecture is strongly connected to the derived invariance of the first Tachikawa conjecture and the Nakayama conjecture. The derived invariance of the Nakayama conjecture was conjectured in \cite[5.2 Conjecture (3)]{MR3551903}.

\section*{Acknowledgements}
The authors would like to thank the anonymous referee for useful comments and suggestions.

The authors would like to thank Hongxing Chen for his comments on an earlier version of this manuscript.

%\enlargethispage{1\baselineskip}

%\interlinepenalty=1000000

\bibliographystyle{alphaurl}
\bibliography{bibarticle}

%\Address

\end{document}